\documentclass{article}

\usepackage{epsfig}
\usepackage{latexsym}

\usepackage{epsfig}
\usepackage{latexsym}
\usepackage{amsmath,amssymb}
\input{epsf}

\newtheorem{prelem}{{\bf Theorem}}

\newtheorem{theorem}{Theorem}[section]
\newtheorem{corollary}[theorem]{Corollary}
\newtheorem{definition}[theorem]{Definition}
\newtheorem{conjecture}[theorem]{Conjecture}

\newtheorem{lemma}[theorem]{Lemma}
\newtheorem{proposition}[theorem]{Proposition}

\newtheorem{observation}[theorem]{Observation}

\newenvironment{proof}{{\bf Proof.}}{\hfill\rule{2mm}{2mm}}
\newtheorem{remarka}[theorem]{Remark}
\newenvironment{remark}{\begin{remarka}\rm}{\hfill\rule{2mm}{2mm}\end{remarka}}
\newtheorem{examplea}[theorem]{Example}
\newenvironment{example}{\begin{examplea}\rm}{\hfill\rule{2mm}{2mm}\end{examplea}}

\def\E {\mathbb{E}}

\def\mod {{\rm mod}}
\def\supp {{\rm supp}}
\def\Re {{\rm Re}}
\def\Tr {{\rm Tr}}

\def\range {{\rm ran}}

\title{On generalizations of Gowers norms and their geometry}
\author{
{\bf  Hamed Hatami} \\
{\small\it Department of Computer Science}\\
{\small University of Toronto} \\
{\small e-mail: hamed@cs.toronto.edu}}
\begin{document}
\maketitle

\tableofcontents

\begin{abstract}
Motivated by the definition of the Gowers uniformity norms, we
introduce and study a wide class of norms. Our aim is to establish
them as a natural generalization of the $L_p$ norms. We shall prove
that these normed spaces share many of the nice properties of the
$L_p$ spaces. Some examples of these norms are $L_p$ norms, trace
norms $S_p$ when $p$ is an even integer, and Gowers uniformity
norms.

Every such norm is defined through a pair of weighted hypergraphs.
In regard to a question of L\'aszl\'o Lov\'asz, we prove several
results in the direction of characterizing all hypergraph pairs that
correspond to norms.
\end{abstract}

\newline
{{\sc Keywords:} Gowers norms, graph densities, hypergraphs, uniform
smoothness, uniform convexity.

\section{Introduction}
 Consider a measurable
function $f:[0,1] \rightarrow \mathbb{C}$. For $1 \le p <\infty$,
the $L_p$ norm of $f$ is defined as
\begin{equation}
\label{eq:Lp} \|f\|_p =  \left(\int |f(x)|^p dx \right)^{1/p} =
\left(\int f(x)^{p/2} \overline{f(x)}^{p/2} dx \right)^{1/p}.
\end{equation}
Next consider a measurable function $f:[0,1]^2 \rightarrow
\mathbb{C}$. The Gowers $2$-uniformity norm of $f$ is defined as
\begin{equation}
\label{eq:U2} \|f\|_{U_2} = \left(\int
f(x_0,y_0)f(x_1,y_1)\overline{f(x_0,y_1)f(x_1,y_0)} dx_0 dx_1 dy_0 dy_1\right)^{1/4}.
\end{equation}
Note that there are similarities between (\ref{eq:Lp}) and
(\ref{eq:U2}): Their underlying vector space is a function space,
and the norm of a function $f$ is defined by a formula of the form
$\left(\int \Pi\right)^{1/p}$, where $p>0$ and $\Pi$ is a product
which involves different copies of powers of $f$ and $\overline{f}$.
The purpose of this article is to use a common framework to study
the norms that are defined in a similar fashion. Our aim is to
establish this class of norms as a natural generalization of the
$L_p$ norms. We shall prove that they share many of the nice
properties of the $L_p$ norms.

An important class of norms that fall into our setting are Gowers
norms. They are introduced by Gowers~\cite{Gowers98,MR2373376} as a
measurement of pseudo-randomness in his proof for  Sz\'{e}meredi's
theorem on arithmetic progressions. The discovery of these norms
resulted in a better understanding of the concept of
pseudo-randomness, and this led to an enormous amount of progress in
the area, and establishment of remarkable results such as Green and
Tao's theorem~\cite{GreenTao} that the primes contain arbitrarily
long arithmetic progressions. Although Gowers norms are very special
case of our framework, surprisingly some of their key properties,
and ideas from pseudo-randomness theory will be needed in our
proofs.

For now let us focus on two-variable functions $f:[0,1]^2
\rightarrow \mathbb{C}$. For finite sets $V_1,V_2$ and functions
$\alpha,\beta:V_1 \times V_2 \rightarrow \mathbb{R}^+$, consider
$$\|f\|_{(\alpha,\beta)} := \left(\int \prod_{(i,j) \in V_1 \times V_2} f(x_i,y_j)^{\alpha(i,j)} \prod_{(i,j) \in V_1 \times V_2} \overline{f(x_i,y_j)}^{\beta(i,j)} \right)^{1/t},$$
where $t:=\sum_{(i,j) \in V} \alpha(i,j)+\beta(i,j)$. A natural
question is that for which $\alpha,\beta$, the function
$\|\cdot\|_{(\alpha,\beta)}$ defines a norm. For example both
formulas
\begin{equation}
\label{eq:introDummy1} \|f\|_{2U_2} :=
\|f^2\|_{U_2}^{1/2}=\left(\int
f(x_0,y_0)^2f(x_1,y_1)^2\overline{f(x_0,y_1)}^2\overline{f(x_1,y_0)}^2
dx_0 dx_1 dy_0 dy_1\right)^{1/8},
\end{equation}
and
\begin{equation}
\label{eq:introDummy2} \left(\int
|f(x_0,y_0)|^{\sqrt{2}}|f(x_1,y_1)|^{\sqrt{2}}
|f(x_0,y_1)||f(x_1,y_0)| dx_0 dx_1 dy_0
dy_1\right)^{1/(2\sqrt{2}+2)},
\end{equation}
can be defined as $\|\cdot\|_{(\alpha,\beta)}$ for proper choices of
functions $\alpha$ and $\beta$. They are both always nonnegative,
and homogenous with respect to scaling. But do they satisfy the
triangle inequality? One of our main results,
Theorem~\ref{thm:outRuling}, says that if $\|\cdot
\|_{(\alpha,\beta)}$ satisfies the triangle inequality, then one of
the following two conditions hold:
\begin{itemize}
\item {\bf Type I:} There exists a constant $s \ge 1$ such that
$\alpha(i,j)=\beta(i,j) \in \{0,s/2\}$, for every $(i,j) \in V_1
\times V_2$;

\item {\bf Type II:} For every $(i,j) \in V_1 \times V_2$,
$\alpha(i,j)=\beta(i,j)=0$, $\alpha(i,j)=1-\beta(i,j)=0$, or
$1-\alpha(i,j)=\beta(i,j)=0$.
\end{itemize}
It follows from the above theorem that neither of
(\ref{eq:introDummy1}) and (\ref{eq:introDummy2}) satisfies the
triangle inequality.  The $L_p$ norm $\|f\|_p = (\int
|f(x,y)|^p)^{1/p}$ is an example of a norm of Type I, and
$\|\cdot\|_{U_2}$ defined in (\ref{eq:U2}) is an example of a norm
of Type II.

Among the key ingredients in the proof of
Theorem~\ref{thm:outRuling} is a H\"older type inequality that we
prove in Lemma~\ref{lem:Holder}. This inequality is extremely useful
in this article and shall be applied frequently. One can think of it
as a common generalization of the classical H\"older inequality and
the Gowers-Cauchy-Schwarz inequality.

We also study the norms $\|\cdot\|_{(\alpha,\beta)}$  from a
geometric point of view, and determine their moduli of smoothness
and convexity. These two parameters are among the most important
invariants in Banach space theory. Our results in particular
determine the moduli of smoothness and convexity of Gowers norms.
They also provide a unified proof for some previously known facts
about $L_p$ and Schatten spaces, and generalize them to a wider
class of norms. When the norm is of Type II we can show that the
corresponding normed space satisfies the so called Hanner
inequality. This inequality has been proven to hold only for a few
spaces, namely the $L_p$ spaces by Hanner~\cite{Hanner}, and the
Schatten spaces $S_p$ for $p \ge 4$ and $1 \le p \le 4/3$ by Ball,
Carlen and Lieb~\cite{BallCarlenLieb}.  We also prove a complex
interpolation theorem for normed spaces of Type I, and use it
together with the Hanner inequality to obtain various optimum
results in terms of the constants involved in the definition of
moduli of smoothness and convexity.

The norms studied  here are generalizations of the graph norms
studied in~\cite{GraphNorms}. For an integer $k>0$, it is well-known
that the $2k$-trace norm of a matrix can be defined through the
graph $C_{2k}$, the cycle of length $2k$. This gives a combinatorial
interpretation of the $2k$-trace norm with many applications in
graph theory. A remarkable recent example is the work of Bourgain
and Gamburd~\cite{BourgainGamburd} on expanders. Inspired by the
fact that the cycles of even length correspond to norms, and the
numerous applications of these norms in graph theory, L\'aszl\'o
Lov\'asz posed the problem of characterizing all graphs that
correspond to norms. The study of this problem is initiated by the
author in~\cite{GraphNorms}, where among other things, a rather
surprising application to Erd\"os-Simonovits-Sidorenko conjecture
has been proven.

Although the framework of the present article is a generalization of
\cite{GraphNorms}, almost all of the results proven here are new
even in the context of the graph norms. In particular we settle an
open question posed in~\cite{GraphNorms}.

\subsection{Notations and Definitions \label{sec:notation}}

In this section we give the formal definition of a hypergraph pair,
and introduce the notations and conventions used throughout the
article. A measure in this article is always a positive measure.

Let $k>0$ be an integer,  $V_1,\ldots,V_k$ be finite nonempty sets
and $V:=V_1 \times \ldots \times V_k$. For $\alpha, \beta:V
\rightarrow \mathbb{R}$, the pair $H=(\alpha,\beta)$ is called a
$k$-\emph{hypergraph pair}.  The size of $H$ is defined as
$$|H|:=\sum_{\omega \in V} |\alpha(\omega)|+|\beta(\omega)|.$$
When we say $H=(\alpha,\beta)$ takes only integer values,  we mean
that $\range(\alpha),\range(\beta) \subseteq \mathbb{Z}$.

 Consider two
$k$-hypergraph pairs: $H=(\alpha,\beta)$ over $V=V_1 \times \ldots
\times V_k$, and $H'=(\alpha',\beta')$ over $W=W_1 \times \ldots
\times W_k$. An \emph{isomorphism} from $H$ to $H'$ is a $k$-tuple
$h=(h_1,\ldots,h_k)$ such that $h_i:V_i \rightarrow W_i$ are
bijections satisfying
$$\alpha(\omega)= \alpha'(h(\omega)), \qquad \beta(\omega)=
\beta'(h(\omega)),$$
for every $\omega=(\omega_1,\ldots,\omega_k) \in V$, where
$h(\omega):= (h_1(\omega_1),\ldots,h_k(\omega_k))$. We say $H$ is
\emph{isomorphic} to $H'$, and denote it by $H \cong H'$, if there
exists an isomorphism from $H$ to $H'$.

Let ${\cal M}=(\Omega,{\cal F},\mu)$ be a measure space.  Every
$\omega \in V$ defines a projection from $\Omega^{V_1} \times \ldots
\times \Omega^{V_k}$ to $\Omega^k$ in a natural way. For a
measurable function $f:\Omega^k \rightarrow \mathbb{C}$, let
$f^H:\Omega^{V_1} \times \ldots \times \Omega^{V_k} \rightarrow
\mathbb{C}$ be defined as
$$f^H(x)  := \left( \prod_{\omega \in V}
f(\omega(x))^{\alpha(\omega)} \right) \left( \prod_{\omega \in V}
\overline{f(\omega(x))^{\beta(\omega)}} \right),$$
where here, and in the sequel we always assume $0^0=1$. As we
discussed above we want to use hypergraph pairs to construct normed
spaces.
\begin{definition}
\label{def:productNorms} Consider a $k$-hypergraph pair
$H=(\alpha,\beta)$ with $\alpha,\beta \ge 0$, and a measure space
${\cal M}=(\Omega,{\cal F},\mu)$. Let $L_H({\cal M})$ be the set of
functions $f:\Omega^k \rightarrow \mathbb{C}$ with $\| \ |f| \
\|_H<\infty$, where for a measurable function $f:\Omega^k
\rightarrow \mathbb{C}$,
\begin{equation}
\label{eq:defProductNorm} \|f\|_H := \left( \int f^H
\right)^{1/|H|}.
\end{equation}
A hypergraph pair is called norming (semi-norming), if $\|\cdot\|_H$
defines a norm (semi-norm) on $L_H({\cal M})$ for every measure
space ${\cal M}=(\Omega,{\cal F},\mu)$.
\end{definition}
\begin{remark}
As the reader might have noticed, the variables and the
infinitesimals are missing from the integral in
(\ref{eq:defProductNorm}). To keep the notation simple, here and in
the sequel when there is no ambiguity we will omit the variables and
infinitesimals from the integrals.
\end{remark}

\begin{remark}
Note that if $H \cong H'$, then for every function $f$ we have $\int
f^H = \int f^{H'}$.
\end{remark}

\begin{remark}
\label{rem:finitization} Note that a hypergraph pair is norming
(semi-norming), if $\|\cdot\|_H$ defines a norm (semi-norm) on
$L_H({\cal M})$, for every measure space ${\cal M}=(\Omega,{\cal
F},\mu)$ with $|\Omega| < \infty$.
\end{remark}

As one would suspect from Definition~\ref{def:productNorms}, the
function $\|\cdot\|_H$ is not a priori a norm. We will pursue the
question: ``Which hypergraph pairs are norming (semi-norming), and
what are the properties of the normed spaces induced by them?''
\begin{remark}
\label{rem:notation} Let $V_1,\ldots,V_k$ be arbitrary finite sets.
For $\psi \in V_1 \times \ldots \times V_k$, we denote by $1_\psi$
the $k$-hypergraph pair $(\delta_\psi,0)$, where $\delta_\psi$ is
the Dirac's delta function: $\delta_{\psi}(\omega)=1$ if
$\omega=\psi$, and $\delta_{\psi}(\omega)=0$ otherwise.

We will apply arithmetic operations to  hypergraph pairs: For
example for two hypergraph pairs $H_1=(\alpha_1,\beta_1)$ and
$H_2=(\alpha_2,\beta_2)$, their sum $H_1+H_2$ and their difference
$H_1-H_2$ are defined respectively as the pairs
$(\alpha_1+\alpha_2,\beta_1+\beta_2)$ and
$(\alpha_1-\alpha_2,\beta_1-\beta_2)$. For a hypergraph pair
$H=(\alpha,\beta)$ define $\overline{H}:=(\beta,\alpha)$, and $r H:=
(r \alpha, r \beta)$ for every $r \in \mathbb{R}$. Now let
$H_1=(\alpha_1,\beta_1)$ be a hypergraph pair over $V_1 \times
\ldots \times V_k$ and $H_2=(\alpha_2,\beta_2)$ be a hypergraph pair
over $W_1 \times \ldots \times W_k$. By considering proper
isomorphisms we can assume that $W_i$ and $V_i$ are all disjoint.
Then the disjoint union $H_1 \dot\cup H_2$ is defined as a
hypergraph pair over $(V_1 \dot\cup W_1) \times \ldots \times (V_k
\dot\cup W_k)$ whose restrictions to $V_1 \times \ldots \times V_k$
and $W_1 \times \ldots \times W_k$ are respectively $H_1$ and $H_2$,
and is defined to be zero everywhere else.  With these definitions,
it is easy to verify that for a measurable function $f:\Omega^k
\rightarrow \mathbb{C}$, we have
\begin{eqnarray*}
f^{H_1+H_2}&=&f^{H_1}f^{H_2} \\
f^{H_1-H_2}&=&f^{H_1}/f^{H_2} \\
f^{\overline{H}}&=&\overline{f^{H}} \\
f^{rH}&=&\left(f^{H}\right)^r=(f^r)^H \\
\int f^{H_1 \dot\cup H_2} &=& \int f^{H_1} \int f^{H_2}.
\end{eqnarray*}
\end{remark}
Consider a hypergraph pair $H$, and note that $\|\cdot\|_{H} = \|
\cdot \|_{H \dot\cup H}$. Thus in order to characterize all norming
(semi-norming) hypergraph pairs it suffices to consider hypergraph
pairs that are minimal according to the following definition:
\begin{definition}
\label{def:minimal} A hypergraph pair $H$ over $V_1 \times \ldots
\times V_k$ is called minimal if
\begin{itemize}
\item For every $i \in [k]$ and $v_i \in V_i$, there exists at least
one $\omega \in \supp(\alpha) \cup \supp(\beta)$ such that
$\omega_i=v_i$.

\item There is no $k$-hypergraph pair $H'$ such that $H \cong H'
\dot\cup H'$.
\end{itemize}
\end{definition}

The next couple of examples show that some well-known families of
normed spaces fall in the framework defined above.
\begin{example}
\label{example:Lp} Let  $L_p=(\alpha,\beta)$ be the $1$-hypergraph
pair defined as $\alpha=\beta=p/2$ over $V_1$ which contains only
one element. Then for a measurable function $f:\Omega \rightarrow
\mathbb{C}$, we have
$$\|f\|_{L_p} = \left(\int f^{p/2} \overline{f^{p/2}}\right)^{1/p}= \left(\int |f|^p\right)^{1/p}=\|f\|_p.$$
Hence in this case the $\|\cdot\|_{L_p}$ norm is the usual $L_{p}$
norm.
\end{example}

\begin{example}
\label{example:Sp} Let $k=2$,  $V_1=V_2=\{0,1,\ldots,m-1\}$, for
some positive integer $m$. Define the $2$-hypergraph pair
$S_{2m}=(\alpha,\beta)$ as
$$\alpha(i,j):=\left\{\begin{array}{lcl}1& \qquad & i=j \\ 0& & \mbox{otherwise} \end{array}\right.$$
$$\beta(i,j):=\left\{\begin{array}{lcl}1& \qquad & i=j+1 (\mod \ m) \\ 0& & \mbox{otherwise} \end{array}\right.$$
Let $\mu$ be the counting measure on a finite set $\Omega$. Then for
$A:\Omega^2 \rightarrow \mathbb{C}$ we have
\begin{eqnarray*}
\|A\|_{S_{2m}}&=& \left(\sum A(x_0,y_0)\overline{A(x_1,y_0)}
A(x_1,y_1)\overline{A(x_2,y_1)} \ldots
A(x_{m-1},y_{m-1})\overline{A(x_0,y_{m-1})}\right)^{1/2m} \\&=&
\left(\Tr(AA^*)^{m}\right)^{1/2m},
\end{eqnarray*}
which shows that in this case the $L_{S_{2m}}$ norm coincides with
the usual $2m$-trace norm of matrices.
\end{example}

\begin{example}
\label{example:Gowers} Let $k$ be a positive integer and
$V_1=\ldots=V_k =\{0,1\}$, and for $\omega \in V_1 \times \ldots
\times V_k$,
$$\alpha(\omega):=\sum_{i=1}^k \omega_i \ (\mod \ 2)$$
and
$$\beta(\omega) := 1 - \alpha(\omega).$$
Then for the $k$-hypergraph pair $U_k=(\alpha,\beta)$,
$\|\cdot\|_{U_k}$ is called the Gowers $k$-uniformity norm.
\end{example}

\subsection{Graph norms and subgraph densities}

Hypergraph norms are important in the study of subgraph and
sub-hypergraph densities. In fact this was one of the main
motivations for studying the graph norms in~\cite{GraphNorms}. We
refer the reader to~\cite{GraphNorms} for the details, but for now
let us define the graph norms  in our notation. Recall that a
bipartite graph is a triple $H=(V_1,V_2,E)$ where $E \subseteq V_1
\times V_2$. Note that every such graph can be identified with a
$2$-hypergraph pair $H=(\alpha,0)$ over $V_1 \times V_2$ where
$\alpha$ is the indicator function of $E$. In~\cite{GraphNorms} two
candidates for being norms are corresponded to $H$. In our notation,
they are defined by the formulas
\begin{equation} \label{eq:graphNorm} \left|
\int f^H \right|^{1/|H|},
\end{equation}
and
\begin{equation} \label{eq:graphAbsoluteNorm} \left(
\int |f|^H \right)^{1/|H|},
\end{equation}
where in (\ref{eq:graphNorm}) $f$ is assumed to be a real-valued
function. In our notation $(\ref{eq:graphNorm})=\|f\|_{H \dot\cup
\overline{H}}$ and $(\ref{eq:graphAbsoluteNorm}) =
\|f\|_{\frac{H+\overline{H}}{2}}$ which shows that our framework in
this article is sufficiently general to include the graph norms.

An important conjecture due to Erd\"os and Simonovits~\cite{Erdos}
(See also Sidorenko~\cite{MR1138091,MR1225933}) can be formulated in
the language of the graph norms. Consider an arbitrary bipartite
graph $H=(V_1,V_2,E)$, a \emph{probability space} ${\cal
P}=(\Omega,{\cal F},\mu)$ and a measurable function $f:\Omega^2
\rightarrow \mathbb{R}^+$. It is conjectured in~\cite{Erdos}
that~\footnote{This form of the conjecture is due to Sidorenko, but
it is equivalent to what is conjectured in~\cite{Erdos}.}
$$\|f\|_1 \le \|f\|_{H}.$$
It has been shown in~\cite{GraphNorms} that if the formula
in~(\ref{eq:graphAbsoluteNorm}) corresponds to a norm, then the
statement of the conjecture is true for $H$. The same arguments hold
in the setting of hypergraph pairs as well, and similar inequalities
can be obtained for norming hypergraph pairs. This follows from
Corollary~\ref{cor:normOrder} below. However it should be noted that
the analogue of Erd\"os-Simonovits-Sidorenko conjecture for
$k$-variable functions where $k>2$ is false (See~\cite{MR1225933}).

The moduli of smoothness and convexity are two dual parameters
assigned to a normed space that play a fundamental role in Banach
space theory. We will discuss them extensively in
Section~\ref{sec:Geometry}. In~\cite{GraphNorms} the moduli of
smoothness and convexity of the normed spaces defined by
(\ref{eq:graphNorm}) are determined, but for the normed spaces
defined by (\ref{eq:graphAbsoluteNorm}) it was left open. This
question will be answered in Theorem~\ref{thm:ConvexitySmoothness}.

\subsection{Constructing norming hypergraph pairs}

The following definition introduces the tensor product of two
hypergraph pairs.

\begin{definition}
Let $H_1=(\alpha_1,\beta_1)$ be a $k$-hypergraph pair over $V_1
\times \ldots \times V_k$ and $H_2=(\alpha_2,\beta_2)$ be a
$k$-hypergraph pair over $W_1 \times \ldots \times W_k$. Then the
tensor product of $H_1$ and $H_2$, is a $k$-hypergraph pair over
$U_1 \times \ldots \times U_k$ where $U_i:=V_i \times W_i$, defined
as $$H_1 \otimes H_2:= (\alpha_1 \otimes \alpha_2 + \beta_1 \otimes
\beta_2, \alpha_1 \otimes \beta_2 + \beta_1 \otimes \alpha_2).$$
\end{definition}

We have already seen in
Examples~\ref{example:Lp},~\ref{example:Sp},~\ref{example:Gowers}
that  norming hypergraph pairs do exist. Theorem~\ref{thm:tensoring}
below shows that it is possible to combine two norming hypergraph
pairs to construct a new one.

\begin{theorem}
\label{thm:tensoring} Let $H_1$ and $H_2$ be two hypergraph pairs.
If $H_1$ and $H_2$ are norming (semi-norming), then  $H_1 \otimes
H_2$ is also norming (semi-norming).
\end{theorem}

The proof of Theorem~\ref{thm:tensoring} is parallel to the proof of
Theorem~2.9 in~\cite{GraphNorms}, and thus we omit it.

The following Lemma which we state without a proof is a
generalization of Theorem 2.8~(ii) in~\cite{GraphNorms}. It can be
easily derived using a similar argument to the proof of
Theorem~2.8~(ii) in~\cite{GraphNorms}.

\begin{lemma}
\label{lem:complete} Consider finite sets $V_1,\ldots,V_k$. For
$\frac{1}{2} \le p < \infty$, the hypergraph pair $K=(p,p)$ over
$V_1 \times \ldots \times V_k$ is norming.
\end{lemma}

\section{Structure of Norming hypergraph pairs \label{sec:Structure}}
In this section we study the structure of  semi-norming hypergraph
pairs. The main result that we  prove in this direction is the
following.

\begin{theorem}
\label{thm:outRuling} Let $H=(\alpha,\beta)$ be a semi-norming
hypergraph pair. Then $H \cong \overline{H}$, and one of the
following two cases hold
\begin{itemize}
\item {\bf Type I}: There exists a real $s \ge 1$, such that for every $\psi \in \supp(\alpha) \cup \supp(\beta)$,
$\alpha(\psi)=\beta(\psi) =s/2$. In this case, $s$ is called the
\emph{parameter} of $H$.

\item {\bf Type II}: For every $\psi \in \supp(\alpha) \cup \supp(\beta)$, we
have $\{\alpha(\psi), \beta(\psi)\}=\{0,1\}$.
\end{itemize}
\end{theorem}

Note that the condition $H \cong \overline{H}$ is trivially
satisfied for every hypergraph pair that satisfies the requirements
of Type~I hypergraph pairs. This is not true for Type~II hypergraph
pairs, and in this case $H \cong \overline{H}$ implies a further
restriction on the structure of the hypergraph pair.

\begin{remark}
Note that if $H$ is of Type I, then for every measure space ${\cal
M}$ and every $f \in L_H({\cal M})$, we have $\| f\|_H = \|\ |f|\
\|_H$. This fact will be used frequently in the sequel.
\end{remark}

Suppose that $H=(\alpha,\beta)$ is a $k$-hypergraph pair over $V_1
\times\ldots \times V_k$. For a subset $S \subseteq [k]$, we use the
notation $\pi_S$ to denote the natural projection from $V_1 \times
\ldots \times V_k$ to $\prod_{i \in S} V_i $. We can construct a
hypergraph pair $H_S:= (\alpha_S, \beta_S)$ where $\alpha_S,\beta_S:
\prod_{i \in S} V_i \rightarrow \mathbb{C}$ are defined as
$$\alpha_S: \omega \mapsto \sum \{\alpha(\omega'): \pi_S(\omega')=\omega\},$$
and
$$\beta_S: \omega \mapsto \sum \{\beta(\omega'): \pi_S(\omega')=\omega\} .$$
By  Remark~\ref{rem:finitization}, we have The following trivial
observation:
\begin{observation}
\label{obs:projection} If  $H=(\alpha,\beta)$  is  a norming
(semi-norming) $k$-hypergraph pair, then for every $S \subseteq
[k]$, $H_S$ is norming (semi-norming).
\end{observation}
\begin{remark}
\label{rem:equalDegrees} The importance of
Observation~\ref{obs:projection} is in that one can apply
Theorem~\ref{thm:outRuling} to $H_S$ to deduce more conditions on
the structure of the original semi-norming hypergraph pair $H$. For
example applying Theorem~\ref{thm:outRuling} to $H_S$ when $S$ has
only one element implies that for every $1 \le i \le k$, there
exists a number $d_i$ such that for every $v_i \in V_i$, we have
$\sum \{\alpha(\omega) : \omega_i=v_i\}=\sum \{\beta(\omega) :
\omega_i=v_i\}=d_i$.
\end{remark}
The next theorem gives another necessary condition on the structure
of a semi-norming hypergraph pair.

\begin{theorem}
\label{thm:spreaded} Suppose that $H=(\alpha,\beta)$ is a
semi-norming $k$-hypergraph pair over $V_1 \times \ldots \times
V_k$. Let $W_i \subseteq V_i$ for $i=1,\ldots,k$, and $H'$ be the
restriction of $H$ to $W_1 \times \ldots \times W_k$. Then
$$\frac{|H'|}{|W_1|+\ldots+|W_k|-1} \le \frac{|H|}{|V_1|+\ldots+|V_k|-1}.$$
\end{theorem}

We present the proofs of Theorems~\ref{thm:outRuling}
and~\ref{thm:spreaded} in Section~\ref{sec:ProofStructure}, but
first we need to develop some tools.

\subsection{Two H\"{o}lder type inequalities}

One of our main tools in the study of hypergraph norms is the trick
of amplification by taking tensor powers. This trick has been used
successfully in many places (see for example~\cite{Ruzsa}).
\begin{definition}
For $f,g:\Omega^k \rightarrow \mathbb{C}$, the tensor product of $f$
and $g$ is defined as $f \otimes g:(\Omega^2)^k \rightarrow
\mathbb{C}$ where $f\otimes g
[(x_1,y_1),\ldots,(x_k,y_k)]=f(x_1,\ldots,x_k)g(y_1,\ldots,y_k)$.
\end{definition}
We have the following trivial observation.
\begin{observation}
\label{obs:tensor} Let $H_1, H_2$  be two $k$-hypergraph pairs, and
$f_1,f_2,g_1,g_2:\Omega^k \rightarrow \mathbb{C}$. Then
$$\int (f_1 \otimes f_2)^{H_1}(g_1 \otimes g_2)^{H_2} = \left(\int f_1^{H_1} g_1^{H_2}\right)\left(\int f_2^{H_1} g_2^{H_2}\right).$$
\end{observation}
Now with Observation~\ref{obs:tensor} in hand, we can prove our
first result about semi-norming hypergraph pairs.
\begin{lemma}
\label{lem:FirstHolder} Let $H=(\alpha,\beta)$ be a semi-norming
hypergraph pair. Then for every measurable space ${\cal M}$, and
every $f,g \in L_H({\cal M})$ the following holds. For every $\psi
\in \supp(\alpha)$,
\begin{equation}
\label{eq:SimpleHolder1} \left|\int f^{H-1_\psi}g^{1_\psi}\right|
\le \|f\|_H^{|H|-1} \|g\|_H,
\end{equation}
and for every $\psi \in \supp(\beta)$
\begin{equation}
\label{eq:SimpleHolder2} \left|\int
f^{H-\overline{1_\psi}}g^{\overline{1_\psi}}\right| \le
\|f\|_H^{|H|-1} \|g\|_H.
\end{equation}
Conversely, if for a measure space ${\cal M}$, and every $f,g \in
L_H({\cal M})$, $\int f^H \in \mathbb{R}^+$, and at least one of
(\ref{eq:SimpleHolder1}) or (\ref{eq:SimpleHolder2}) holds for some
$\psi \in V_1 \times \ldots \times V_k$, then $\|\cdot \|_H$ is a
semi-norm on $L_H({\cal M})$.
\end{lemma}
\begin{proof}
First we prove the converse direction which is easier. Consider two
measurable functions $f,g:\Omega^k \rightarrow \mathbb{C}$ and
suppose that (\ref{eq:SimpleHolder1}) holds for some $\psi \in V_1
\times \ldots \times V_k$. Then
\begin{eqnarray*}
\|f+g\|_H^{|H|} &=& \int (f+g)^H = \int
(f+g)^{H-1_\psi}(f+g)^{1_\psi}\\&=& \int (f+g)^{H-1_\psi}f^{1_\psi}+
\int (f+g)^{H-1_\psi}g^{1_\psi} \le \|f+g\|_H^{|H|-1}\|f\|_H +
\|f+g\|_H^{|H|-1}\|g\|_H,
\end{eqnarray*}
which simplifies to the triangle inequality. The proof of the case
where (\ref{eq:SimpleHolder2}) holds is similar.

Now let us turn to the other direction.  Suppose that $H$ is a
semi-norming hypergraph pair. Consider $f,g \in L_H({\cal M})$. We
might assume that $\|f\|_H \neq 0$, as otherwise one can instead
consider a small perturbation of $f$. Since $\|\cdot\|_H$ is a
semi-norm, for every $t \in \mathbb{R}^+$ and every $f,g:\Omega
\rightarrow \mathbb{C}$, we have $\|f+tg\|_H \le \|f\|_H + t\|g\|_H$
which implies that
\begin{equation}
\label{eq:1derivative} \left. \frac{d \|f+gt\|_H}{dt}\right|_0 \le
\|g\|_H.
\end{equation}
Computing the derivative
$$\frac{d (f+tg)^H}{dt}  =
\sum_{\psi \in \supp(\alpha)} \alpha(\psi)(f+tg)^{H-1_\psi}
g^{1_\psi}+\sum_{\psi \in \supp(\beta)}
\beta(\psi)(f+tg)^{H-\overline{1_\psi}} g^{\overline{1_\psi}},$$
shows that
$$\frac{d \|f+tg\|_H}{dt}  =
\frac{1}{|H|}\|f+tg\|_H^{1-|H|}\left(\int \sum_{\psi \in
\supp(\alpha)} \alpha(\psi)(f+tg)^{H-1_\psi} g^{1_\psi}+\sum_{\psi
\in \supp(\beta)}\beta(\psi)(f+tg)^{H-\overline{1_\psi}}
g^{\overline{1_\psi}}\right).$$
Thus by (\ref{eq:1derivative}),
$$ \frac{1}{|H|}\|f\|_H^{1-|H|}\left(\int \sum_{\psi \in \supp(\alpha)}
\alpha(\psi)f^{H-1_\psi} g^{1_\psi}+ \sum_{\psi \in
\supp(\beta)}\beta(\psi)f^{H-\overline{1_\psi}}
g^{\overline{1_\psi}}\right) \le \|g\|_H,$$
or equivalently
\begin{equation}
\label{eq:preTensor} \frac{1}{|H|}\left(\int \sum_{\psi \in
\supp(\alpha)} \alpha(\psi)f^{H-1_\psi} g^{1_\psi}+ \sum_{\psi \in
\supp(\beta)}\beta(\psi)f^{H-\overline{1_\psi}}
g^{\overline{1_\psi}}\right) \le \|f\|_H^{|H|-1} \|g\|_H.
\end{equation}
Since (\ref{eq:preTensor}) holds for every measure space and every
pair of measurable functions, for every integer $m>0$, we can
replace $f$ and $g$ in~(\ref{eq:preTensor}), respectively with
$f^{\otimes m} \otimes \overline{f}^{\otimes m}$ and $g^{\otimes m}
\otimes \overline{g}^{\otimes m}$, and apply
Observation~\ref{obs:tensor} to obtain
\begin{equation}
\label{eq:postTensor} \frac{1}{|H|}\left(\sum_{\psi \in
\supp(\alpha)} \alpha(\psi) \left|\int f^{H-1_\psi}
g^{1_\psi}\right|^{2m}+ \sum_{\psi \in \supp(\beta)}\beta(\psi)
\left|\int f^{H-\overline{1_\psi}}
g^{\overline{1_\psi}}\right|^{2m}\right) \le \left(\|f\|_H^{|H|-1}
\|g\|_H \right)^{2m}.
\end{equation}
But since~(\ref{eq:postTensor}) holds for every $m$, it establishes
(\ref{eq:SimpleHolder1}) and (\ref{eq:SimpleHolder2}) as
$$\frac{1}{|H|}\left(\sum_{\psi \in \supp(\alpha)} \alpha(\psi) +
\sum_{\psi \in \supp(\alpha)} \alpha(\psi)\right) =1.$$
\end{proof}

We have the following corollary to Lemma~\ref{lem:FirstHolder}.
\begin{corollary}
\label{cor:MoreThanOne} If $H$ is a semi-norming hypergraph pair,
then $\alpha(\omega)+\beta(\omega) \ge 1$, for every $\omega \in
\supp(\alpha) \cup \supp(\beta)$.
\end{corollary}
\begin{proof}
Let the underlying measure space be the set $\{0,1\}$ with the
counting measure. Consider $\omega \in \supp(\alpha)$, and note that
by (\ref{eq:SimpleHolder1}), for every pair of functions
$f,g:\{0,1\}^k \rightarrow \mathbb{C}$, we have
\begin{equation}
\label{eq:CorHolder} \left|\int f^{H-1_\psi}g^{1_\psi}\right| \le
\|f\|_H^{|H|-1} \|g\|_H.
\end{equation}
For every $x=(x_1,\ldots,x_k) \in \{0,1\}^k$, define $g(x):=1$ and
$$f(x):=\left\{
\begin{array}{lcl}
\epsilon& \qquad& x_1=\ldots=x_k=1 \\ 1 &\qquad &\mbox{otherwise}
\end{array} \right.$$
Then $\left|\int f^{H-1_\psi}g^{1_\psi}\right|=\left|\int
f^{H-1_\psi} \right| \ge \epsilon^{\alpha(\omega)+\beta(\omega)-1},$
while $\|f\|_H \le \|g\|_H= \|1\|_H$, which contradicts
(\ref{eq:CorHolder}) for sufficiently small $\epsilon>0$, if
$\alpha(\omega)+\beta(\omega) <1$.
\end{proof}

Under some extra conditions it is possible to extend
(\ref{eq:SimpleHolder1}) and (\ref{eq:SimpleHolder2}) to  a much
more powerful inequality.
\begin{lemma}
\label{lem:Holder} Let $H$ be a semi-norming hypergraph pair, and
$H_1,\ldots,H_n$ be nonzero~\footnote{i.e. $H_i \neq (0,0)$ for
every $1 \le i \le n$.} hypergraph pairs satisfying
$H_1+H_2+\ldots+H_n=H$. Then for every measure space ${\cal M}$ and
functions $f_1,f_2,\ldots,f_n \in L_H({\cal M})$, we have
$$\left| \int f_1^{H_1} f_2^{H_2} \ldots f_n^{H_n} \right| \le \|f_1\|_H^{|H_1|} \ldots \|f_n\|_H^{|H_n|},$$
provided that at least one of the following two conditions hold:
\begin{enumerate}
\item[(a)] We have $f_1,\ldots,f_n \ge 0$.

\item[(b)] For every $H_i=(\alpha_i,\beta_i)$, the functions $\alpha_i,\beta_i$ take only
integer values.
\end{enumerate}
\end{lemma}
\begin{proof}
Let us first assume that $f_1,\ldots,f_n \ge 0$. Suppose to the
contrary that
\begin{equation}
\label{eq:contrary} \int f_1^{H_1} f_2^{H_2} \ldots f_n^{H_n}  >
\|f_1\|_{H}^{|H_1|} \ldots  \|f_n\|_{H}^{|H_n|}.
\end{equation}
After normalization we can assume that
$\|f_1\|_H,\|f_2\|_H,\ldots,\|f_n\|_H \le 1$ while the right-hand
side of (\ref{eq:contrary}) is strictly greater than $1$. Since
(\ref{eq:contrary}) remains valid after small perturbations of
$f_i$'s, without loss of generality we might also assume that for
every $1 \le i \le n$, $f_i$ does not take the zero value on any
point. Consider a positive integer $m$, and note that by
Observation~\ref{obs:tensor}
\begin{eqnarray*}
\int \left(\sum_{i=1}^n f_i^{\otimes m}\right)^H &=& \int
\prod_{i=1}^n \left(\sum_{i=1}^n f_i^{\otimes m}\right)^{H_i} = \int
(f_1^{\otimes m})^{H_1} \ldots (f_n^{\otimes
m})^{H_n}\left(\prod_{i=1}^n\left(\frac{f_1^{\otimes m}+\ldots+f_n^{\otimes m}}{f_i^{\otimes m}}\right)^{H_i}\right)\\
&\ge& \int (f_1^{\otimes m})^{H_1} \ldots (f_n^{\otimes m})^{H_n} =
\left(\int f_1^{H_1} \ldots f_n^{H_n}\right)^m.
\end{eqnarray*}
On the other hand, Observation~\ref{obs:tensor} shows that
$\|f_i^{\otimes m}\|_H = \|f_i\|_H^m \le 1 $ for every $i \in [n]$.
Then for sufficiently large $m$ we get a contradiction:
$$\left\|\sum f_i^{\otimes m}\right\|_H \ge \left(\int f_1^{H_1} \ldots f_n^{H_n}\right)^{m/|H|} > n.$$

Next consider the case where $f_i$ are not necessarily positive, but
we know that $\alpha_i,\beta_i$ all take only integer values. Again
to get a contradiction assume that
$$\left|\int f_1^{H_1} \ldots f_n^{H_n} \right|  > 1 \ge \|f_1\|_{H}^{|H_1|} \ldots \|f_n\|_{H}^{|H_n|}.$$
where $\|f_1\|_{H},\ldots,\|f_n\|_{H} \le 1$. In this case for every
$i \in [n]$, we will consider $f_i^{\otimes m} \otimes
\overline{f_i}^{\otimes m}$. Let ${\cal H}$ denote the set of all
$n$-tuples of nonzero hypergraph pairs $(H'_1,H'_2,\ldots,H'_n)$
where $H'_i$'s take only nonnegative integer values and
$H'_1+H'_2+\ldots+H'_n=H$. By Observation~\ref{obs:tensor}
$$\int \prod_{i=1}^n (f_i^{\otimes m}\otimes \overline{f_i}^{\otimes m})^{H'_i}
= \left|\int \prod_{i=1}^n f_i^{H'_i} \right|^{\otimes 2m} \ge 0.$$
Now by expanding the product defined by $H$, we have
\begin{eqnarray*}
\int \left(\sum_{i=1}^n f_i^{\otimes m}\otimes
\overline{f_i}^{\otimes m}\right)^H &=& \sum_{(H'_1,\ldots,H'_n) \in
{\cal H}}\int \prod_{i=1}^n \left(f_i^{\otimes m}\otimes
\overline{f_i}^{\otimes m}\right)^{H'_i}\\&\ge&\int \prod_{i=1}^n
\left(f_i^{\otimes m}\otimes \overline{f_i}^{\otimes
m}\right)^{H_i}=\left|\int \prod_{i=1}^n f_i^{H_i} \right|^{2m},
\end{eqnarray*}
which leads to a contradiction similar to the previous case.
\end{proof}

\begin{remark}
It is possible to show that Lemma~\ref{lem:Holder} does not
necessarily hold in the general case where none of the two
conditions are satisfied. To see this consider $S_4$ from
Example~\ref{example:Sp}. By Lemma~\ref{lem:FirstHolder}, if
Lemma~\ref{lem:Holder} holds for the decomposition
$S_4=\frac{1}{3}S_4 + \frac{1}{3}S_4 + \frac{1}{3}S_4$, then $3S_4$
would be a semi-norming hypergraph pair. But
Theorem~\ref{thm:outRuling} implies that $3S_4$ is not a
semi-norming hypergraph pair.
\end{remark}

Consider a \emph{probability space} ${\cal P}=(\Omega,{\cal
F},\mu)$. It is well-known that for every $1 \le p \le q $, and for
every $f \in L_q({\cal P})$, we have $\|f\|_p \le \|f\|_q$. The next
corollary generalizes this to hypergraph pairs.

\begin{corollary}
\label{cor:normOrder} Let $H=(\alpha,\beta)$ be a semi-norming
$k$-hypergraph pair. Consider a probability space ${\cal
P}=(\Omega,{\cal F},\mu)$ and $f \in L_H({\cal P})$. Let
$K=(\alpha', \beta')$ be a nonzero $k$-hypergraph pair over the same
domain as $H$ such that $\alpha' \le \alpha$ and $\beta' \le \beta$.
Then
$$\left| \|f\|_{K} \right|\le \|f\|_H,$$
provided that at least one of the  following three conditions holds:
\begin{enumerate}
\item[(a)] $f \ge 0$.
\item[(b)] $H$ is of type I.
\item[(c)] The functions $\alpha,\beta,\alpha',\beta'$ take only integer values.
\end{enumerate}
\end{corollary}
\begin{proof}
Parts (a) and (b) follow from applying Lemma~\ref{lem:Holder}~(a),
with parameters $n:=2$, $H_1:=K$, $H_2:=H-K$, $f_1:=|f|$ and
$f_2:=1$.

Part (c) follows from applying Lemma~\ref{lem:Holder}~(b), with
parameters $n:=2$, $H_1:=K$, $H_2:=H-K$, $f_1:=f$ and $f_2:=1$.
\end{proof}

\subsection{Factorizable hypergraph pairs}
In this section we characterize all norming and semi-norming
$1$-hypergraph pairs. As it is mentioned before, it suffices to
consider the hypergraph pairs that are minimal according to
Definition~\ref{def:minimal}. We have already seen one class of
examples of norming $1$-hypergraph pairs, namely the $1$-hypergraph
pairs $L_p$ of Example~\ref{example:Lp}. There exists also a
semi-norming $1$-hypergraph pair that is not norming. Let $G=(1,0)$
be the $1$-hypergraph pair over a set $V_1$ of size $1$. Then for a
measure space ${\cal M}=(\Omega,{\cal F},\mu)$ and a measurable
$f:\Omega \rightarrow \mathbb{C}$ we have $\|f\|_{G \dot\cup
\overline{G}} = \left|\int f \right|$ which defines a semi-norm. The
next proposition shows that  these are the only examples.

\begin{proposition}
\label{prop:1products} If $H$ is a \emph{minimal} norming
$1$-hypergraph pair, then there exists $1 \le p < \infty$ such that
$H \cong L_p$. If $H$ is a minimal semi-norming $1$-hypergraph pair
that is not norming, then  $H \cong G \dot\cup \overline{G}$, where
$G=(1,0)$ is a $1$-hypergraph pair over a set $V_1$ of size $1$.
\end{proposition}
To prove Proposition~\ref{prop:1products} we need to study the
hypergraph pairs which are decomposable into disjoint union of other
hypergraph pairs.
\begin{definition}
A hypergraph pair $H=(\alpha,\beta)$ is called \emph{factorizable},
if it is the disjoint union of two hypergraph pairs.
\end{definition}

The next proposition shows that two non-factorizable hypergraph
pairs define identical norms, if and only if they are isomorphic.
For the proof, we need an easy fact stated in the following Remark.

\begin{remark}
\label{rem:polynomials} Let $x_1,\ldots,x_n$ be $n$ complex
variables. Define a term as a product $\prod_{i=1}^n x_i^{p_i}
\overline{x_i}^{q_i}$, where $p_i,q_i$ are nonnegative reals. Now
let $P$ and $Q$ be two formal finite sums of terms. It is easy to
see that $P$ and $Q$ are equal as functions on $\mathbb{C}^n$, if
and only if they are equal as formal sums.
\end{remark}

\begin{proposition}
\label{prop:isometry} Let $H_1$ and $H_2$ be two minimal
$k$-hypergraph pairs. Suppose that either $H_1$ and $H_2$ are both
non-factorizable, or we have $|H_1|=|H_2|$. Then
\begin{itemize}
\item If for every measure space $(\Omega,{\cal F},\mu)$, and every $f:\Omega^k \rightarrow \mathbb{C}$,
      $\|f\|_{H_1} =\|f\|_{H_2}$, then $H_1 \cong H_2$.

\item If for every measure space $(\Omega,{\cal F},\mu)$, and every $f:\Omega^k \rightarrow \mathbb{C}$,
      $ \|f\|_{H_1} =\overline{ \|f\|_{H_2}}$, then $H_1 \cong \overline{H_2}$.
\end{itemize}
\end{proposition}
\begin{proof}
Suppose that $H_1$ and $H_2$ are respectively defined over $V_1
\times \ldots \times V_k$ and $W_1 \times \ldots \times W_k$.  First
assume that $H_1$ and $H_2$ are both non-factorizable. Let $\mu$ be
the counting measure on $\Omega=[m]$, where $m> \sum_{i=1}^k |V_i| +
|W_i|$ is a positive integer. Suppose that for every $f:\Omega^k
\rightarrow \mathbb{C}$ with have $\|f\|_{H_1}=\|f\|_{H_2}$. Then
define $f(x_1,\ldots,x_k)$ to be equal to $1$, if $x_1=\ldots=x_k$,
and equal to $0$ otherwise. Since $H_1$ and $H_2$ are
non-factorizable it is easy to see that $\int f^{H_1}=\int
f^{H_2}=|\Omega|$ and we deduce that $|H_1|=|H_2|$. So it is
sufficient to prove the proposition for the case where
$|H_1|=|H_2|$.

Now for every $f:\Omega^k \rightarrow \mathbb{C}$, we have $\int
f^{H_1} = \int f^{H_2}$. For $1 \le i \le k$, consider $f_i:\Omega^k
\rightarrow \{0,1\}$ defined as $f_i(x_1,\ldots,x_k)=1$ if and only
if $x_1=\ldots=x_{i-1}=x_{i+1}=\ldots=x_{k}=1$. Then it is easy to
see that $\int f_i^{H_1}=|\Omega|^{|V_i|}$ and $\int
f_i^{H_2}=|\Omega|^{|W_i|}$ which implies $|V_i|=|W_i|$. Thus
without loss of generality we may assume that
$V_i=W_i=\{1,\ldots,|V_i|\}$, for every $1 \le i \le k$. Now for
every $f:\Omega^k \rightarrow \mathbb{C}$ we have
\begin{equation}
\label{eq:polyEquality} \sum_{x \in \Omega^{V_1} \times \ldots
\times \Omega^{V_k}} \prod_{\omega \in V}
f(\omega(x))^{\alpha(\omega)}
\overline{f(\omega(x))}^{\beta(\omega)}=\sum_{x \in \Omega^{V_1}
\times \ldots \times \Omega^{V_k}} \prod_{\omega \in V}
f(\omega(x))^{\alpha'(\omega)}
\overline{f(\omega(x))}^{\beta'(\omega)}
\end{equation}
Consider
$x=[(1,\ldots,|V_1|),(1,\ldots,|V_2|),\ldots,(1,\ldots,|V_k|)] \in
\Omega^{V_1} \times \ldots \times \Omega^{V_k}$. Then
$\omega(x)=\omega$ for every $\omega \in V$, and hence
\begin{equation}
\label{eq:polyEqualityX} \prod_{\omega \in V}
f(\omega(x))^{\alpha(\omega)}
\overline{f(\omega(x))}^{\beta(\omega)} = \prod_{\omega \in V}
f(\omega)^{\alpha(\omega)} \overline{f(\omega)}^{\beta(\omega)}.
\end{equation}
Since (\ref{eq:polyEqualityX}) appears in the sum in the left-hand
side of (\ref{eq:polyEquality}), by Remark~\ref{rem:polynomials} it
must also appear as a term in the right-hand side of
(\ref{eq:polyEquality}). Hence there exists
$y=[(y_{1,1},\ldots,y_{1,|V_1|}),\ldots,(y_{k,1},\ldots,y_{k,|V_k|})]
\in \Omega^{V_1} \times \ldots \times \Omega^{V_k}$ such that
\begin{equation}
\prod_{\omega \in V} f(\omega(y))^{\alpha'(\omega)}
\overline{f(\omega(y))}^{\beta'(\omega)} = \prod_{\omega \in V}
f(\omega)^{\alpha(\omega)} \overline{f(\omega)}^{\beta(\omega)}.
\end{equation}
By minimality (see Definition~\ref{def:minimal}), for every $v \in
V_i$, there exists $\omega=(\omega_1,\ldots,\omega_k) \in
\supp(\alpha) \cup \supp(\beta)$ such that $\omega_i=v$. This
implies $\{y_{i,1},\ldots,y_{i,|V_i|}\}=V_i$, for every $1 \le i \le
k$. Now $h=(h_1,\ldots,h_k)$ defined as $h_i : j \mapsto y_{i,j}$
(for every $1 \le i \le k$ and $1 \le j \le |V_i|$) is an
isomorphism between $H_1$ and $H_2$.

In the second part of the proposition where it is assumed
$\|f\|_{H_1} = \overline{\|f\|_{H_2}}$, instead
of~(\ref{eq:polyEquality}) one obtains that the left-hand side
of~(\ref{eq:polyEquality}) is equal to the conjugate of the
right-hand side. The proof then proceeds similar to the previous
case.
\end{proof}

\begin{theorem}
\label{thm:factor} Let $H=H_1 \dot\cup H_2 \dot\cup \ldots \dot\cup
H_m$ be a semi-norming hypergraph pair such that $H_i$ are all
non-factorizable. Then for every measure space ${\cal M}$ and every
$f \in L_H({\cal M})$ we have $$\|f\|_{H_1 \dot\cup
\overline{H_1}}=\|f\|_{H_2 \dot\cup \overline{H_2}}= \ldots =
\|f\|_{H_m \dot\cup \overline{H_m}} =\|f\|_H.$$
\end{theorem}
\begin{proof}
Let $H=G_1 \dot\cup G_2$ be semi-norming, where $G_1$ and $G_2$ are
not necessarily non-factorizable, ${\cal M}=(\Omega,{\cal F},\mu)$
be a measure space, and $f \in L_H({\cal M})$. Note that
$$ \|f\|_H = \|f\|_{G_1}^{\frac{|G_1|}{|H|}} \|f\|_{G_2}^{\frac{|G_2|}{|H|}}=
\|f\|_{G_1}^{\frac{|G_1|}{|G_1|+|G_2|}}
\|f\|_{G_2}^{\frac{|G_2|}{|G_1|+|G_2|}}.$$
It follows from Theorem~\ref{thm:outRuling} that either $H$ is of
Type I, or $H$ and $G_1$ both take only integer values. Hence by
Corollary~\ref{cor:normOrder}
$$\left| \|f\|_{G_1}^{|G_1|} \right| \le
\|f\|_H^{|G_1|}=\|f\|_{G_1}^{\frac{|G_1|^2}{|G_1|+|G_2|}}\|f\|_{G_2}^{\frac{|G_2|^2}{|G_1|+|G_2|}},$$
which simplifies to
$$\left| \|f\|_{G_1} \right|\le \left| \|f\|_{G_2} \right|.$$
Similarly one can show that $\left| \|f\|_{G_2} \right| \le \left|
\|f\|_{G_1} \right|$ and thus $\left| \|f\|_{G_1} \right|= \left|
\|f\|_{G_2} \right|$.

By induction we conclude that $|\|f\|_{H_1}|=\ldots =|\|f\|_{H_2}|$,
for every measure space ${\cal M}=(\Omega,{\cal F},\mu)$ and every
$f \in L_H({\cal M})$, and this completes the proof as $\left|
\|f\|_{H_i} \right| = \|f\|_{H_i \dot\cup \overline{H_i}}$.
\end{proof}

Now we can state the proof of Proposition~\ref{prop:1products}.

\noindent \begin{proof}[Proposition~\ref{prop:1products}] Consider a
semi-norming $1$-hypergraph pair $H$ over a set
$V_1=\{v_1,\ldots,v_m\}$. Consider the factorization $H=H_1 \dot\cup
H_2 \dot\cup \ldots \dot\cup H_m$, where $H_i$ is a $1$-hypergraph
pair over $\{v_i\}$. By Theorem~\ref{thm:factor}, always $\|f\|_{H_1
\dot\cup \overline{H_1}}=\|f\|_{H_2 \dot\cup \overline{H_2}}= \ldots
= \|f\|_{H_m \dot\cup \overline{H_m}} =\|f\|_H$. By
Theorem~\ref{thm:outRuling}, for every $1 \le i \le m$, either $H_i
\dot\cup \overline{H_i} \cong L_p \dot\cup \overline{L_p}$ for some
$1 \le p <\infty$, or $H_i \dot\cup \overline{H_i} \cong G \dot\cup
\overline{G}$ which completes the proof.
\end{proof}

\subsection{Semi-norming hypergraph pairs that are not norming}
In this section we study the structure of the semi-norming
hypergraph pairs which are \emph{not} norming. Consider a
semi-norming $k$-hypergraph pair $H=(\alpha,\beta)$ over $V:=V_1
\times \ldots \times V_k$ of Type I with parameter $s=2m$, where $m$
is a positive integer. Since $H$ is of Type I, it is trivially
norming. Consider an arbitrary positive integer $k'$. We want to use
$H$ to construct a semi-norming $(k+k')$-hypergraph pair that is not
norming. For $k+1 \le i \le k+k'$, let $V_i:= \supp(\alpha) \times
\{1,\ldots,s\}$. Now $G=(\alpha',\beta')$ is defined by
$$
\alpha(v_1,\ldots,v_{k+k'}):=\left\{\begin{array}{lcl}1 &\quad &
\mbox{$v_{k+1}=\ldots=v_{k+k'}=([v_1,\ldots,v_k],i)$ where $1 \le i
\le m$} \\
0 & \qquad & \mbox{otherwise}
\end{array} \right.
$$
and
$$
\beta(v_1,\ldots,v_{k+k'}):=\left\{\begin{array}{lcl}1 &\quad &
\mbox{$v_{k+1}=\ldots=v_{k+k'}=([v_1,\ldots,v_k],i)$ where $m+1 \le
i \le 2m$} \\
0 & \qquad & \mbox{otherwise}
\end{array} \right.
$$
Consider a measure space ${\cal M}=(\Omega,{\cal F},\mu)$, and an
integrable  function $f:\Omega^{k+k'} \rightarrow \mathbb{C}$. Let
$F:\Omega^k \rightarrow \mathbb{C}$ be defined as $F(x_1,\ldots,x_k)
= \int f(x_1,\ldots,x_{k+k'}) dx_{k+1} \ldots dx_{k+k'}$. It is not
difficult to see that $\|f\|_G = \|F\|_H$, which shows that $G$ is
semi-norming. On the other-hand if $\int f dx_{k+1} \ldots
dx_{k+k'}=0$, then $\|f\|_G=\|F\|_H=\|0\|_H=0$ which implies that
$G$ is not norming. The next proposition shows that in fact every
semi-norming hypergraph pair which is not norming is of this form.

\begin{proposition}
Let $H=(\alpha,\beta)$ be a semi-norming $k$-hypergraph pair of Type
II over $V:=V_1 \times \ldots \times V_k$. Define $S$ to be the set
of all $ 1 \le i \le k$ such that for every $v_i \in V_i$,
$$\sum \{ \alpha(\omega)+\beta(\omega) : \omega \in V,
\omega_i=v_i\}=1.$$
Then $H_{[k] \setminus S}$ is a norming hypergraph pair of Type I.
\end{proposition}
\begin{proof}
Consider a measure space ${\cal M}=(\Omega, {\cal F},\mu)$. Note
that if $S \neq \emptyset$, then for every $i \in S$, every $f \in
L_H({\cal M})$ with $\int f(x_1,\ldots,x_k) dx_i=0$ satisfies
$\|f\|_H=0$. So if $H$ is norming, then $H_{[k]\setminus S}=H$, and
the proposition holds. Consider a $k$-hypergraph pair
$H=(\alpha,\beta)$ over $V:=V_1 \times \ldots \times V_k$ which is
not norming. Then there exists a function $f \in L_H({\cal M})$, for
some measure space ${\cal M}=(\Omega,{\cal F},\mu)$, such that $\int
f^H=0$ and $f \neq 0$. Lemma~\ref{lem:FirstHolder}, then shows that
for every $g \in L_H({\cal M})$, and every $\psi \in \supp(\alpha)$,
\begin{equation}
\label{eq:zeroProd} \int g^{H - 1_\psi} f^{1_\psi} =0.
\end{equation}
Since $f \neq 0$, there exists measurable sets
$\Gamma_1,\ldots,\Gamma_k \subseteq \Omega$ such that
$\int_{\Gamma_1 \times \ldots \times \Gamma_k} f \neq 0$. Define
$g:\Omega^k \rightarrow \{0,1\}$, as
$$g(x_1,\ldots,x_k) =
\left\{\begin{array}{lcl}
1 & \qquad & (x_1,\ldots,x_k) \in \Gamma_1 \times \ldots \times \Gamma_k \\
0 & &\mbox{otherwise}
\end{array} \right. $$
Suppose that for every $i \in [k]$, there exists $\omega \in
\supp(\alpha) \cup \supp(\beta)$ such that $\omega \neq \psi$ but
$\omega_i=\psi_i$. Then it is easy to see that for every $x \in
\Omega^{V_1} \times \ldots \times \Omega^{V_k}$,
$$g^{H-1_\psi}(x)=\left\{\begin{array}{lcl} 1 &\qquad  &\psi(x) \in \Gamma_1 \times \ldots \times
\Gamma_k \\ 0 & &{\rm otherwise}\end{array}\right.$$
 But then $\int g^{H - 1_\psi} f^{1_\psi} = \int_{\Gamma_1
\times \ldots \times \Gamma_k} f  \neq 0$ contradicting
(\ref{eq:zeroProd}).

 It follows from (\ref{eq:zeroProd}) and its analogue
for $\psi \in \supp(\beta)$ that the following holds: For every
$\psi=(\psi_1,\ldots,\psi_k) \in \supp(\alpha) \cup \supp(\beta)$,
there exists $i \in [k]$ such that
$$\{ \omega \in \supp(\alpha)
\cup \supp(\beta): \mbox{$\omega_i=\psi_i$ and $\omega \neq \psi$}
\}= \emptyset,$$
or in other words: $\sum\{\alpha(\omega)+\beta(\omega): \omega \in
V, \omega_i=\psi_i\}=1$.  Now Remark~\ref{rem:equalDegrees} shows
that $i \in S$. By Observation~\ref{obs:projection}, $H_{[k]
\setminus S}$ is semi-norming, but then  maximality of $S$  shows
that it is also norming.
\end{proof}
\subsection{Some facts about Gowers norms}
In this section we prove some facts about Gowers norms that are
needed in the subsequent sections. These facts are only proved as
auxiliary results, and thus our aim is not to obtain the best
possible bounds or to prove them in the most general possible
setting.

Let $V_1=\ldots=V_k=\{0,1\}$, and $U_k$ be the Gowers $k$-hypergraph
pair defined in Example~\ref{example:Gowers}. Consider a measure
space ${\cal M}=(\Omega,{\cal F},\mu)$ and measurable functions
$f_\omega: \Omega^k \rightarrow \mathbb{C}$ for $\omega \in V:=V_1
\times \ldots \times V_k$. The following inequality due to
Gowers~\cite{Gowers98} (see also~\cite{TaoMontreal}) can be proven
by iterated applications of the Cauchy-Schwarz inequality:
\begin{equation}
\left| \int \prod_{\omega \in V} f_\omega^{1_\omega} \right|\le
\prod_{\omega \in V} \|f_\omega\|_{U_k}.
\end{equation}
Since always $\|f\|_{U_k} \le \|f\|_\infty$, we have the following
easy corollary.
\begin{corollary}
\label{cor:Gowers} Let $H=(\alpha,\beta)$ be a $k$-hypergraph pair
over $W:=W_1 \times W_2 \times \ldots \times W_k$, and $\psi \in W$
be such that $\alpha(\psi) = \beta(\psi) =0$. Then for the measure
space ${\cal M}=(\Omega,{\cal F},\mu)$ and every pair of measurable
functions $f,g:\Omega^k \rightarrow \mathbb{C}$, we have
$$\left|\int f^H g^{1_\psi} \right| \le \|g\|_{U_k} \|f\|_\infty^{|H|}.$$
\end{corollary}

The next Lemma shows that there exists a function $g$ such that its
range is $\{-1,1\}$ but its Gowers norm is arbitrarily small.

\begin{lemma}
\label{lem:RandomFunc} For every $\epsilon>0$, there exists a
\emph{probability} space $(\Omega,{\cal F},\mu)$ and a function
$g:\Omega^k \rightarrow \{-1,1\}$ such that $\|g\|_{U_k} \le
\epsilon$  and $\int g=0$.
\end{lemma}
\begin{proof}
Consider a sufficiently large even integer $m$, set $\Omega=[m]$,
and let $\mu$ be the uniform probability measure on $\Omega$. Define
$g$ randomly so that $\{g(\omega)\}_{\omega \in \Omega^k}$ are
independent Bernoulli random variables taking values uniformly
 in $\{-1,1\}$. Then it is easy to see that
$$ \mbox{$\E \left(\int g\right)^2=o_{m \rightarrow \infty}(1)$ \qquad and \qquad $\E \left(\int g^{U_k}\right)^2 = o_{m \rightarrow \infty}(1)$}.$$
Hence for sufficiently large $m$, there exists $g_0:\Omega^k
\rightarrow \{-1,1\}$ such that $|\int g_0| \le (\epsilon/4)^{2^k}$
and $\|g_0\|_{U_k} \le \epsilon/2$. Trivially there exists
$g_1:\Omega^k \rightarrow \{-1,1\}$ such that $\int g_1=0$ and $\int
|g_1-g_0| \le (\epsilon/4)^{2^k}$. Then by H\"older's inequality
$$\|g_0-g_1\|_{U_k} = \left(\int (g_0-g_1)^{U_k}\right)^{-2^k} \le \|g_0-g_1\|_{2^k}  \le 2 (\epsilon/4)=\epsilon/2,$$
where in the last inequality we used the fact that the range of
$g_0-g_1$ is $\{-2,0,2\}$. Now
$$\|g_1\|_{U_k} \le \|g_0\|_{U_k} + \|g_0-g_1\|_{U_k} \le \epsilon,$$
which shows that $g_1$ is the desired function.
\end{proof}

\begin{lemma}
\label{lem:subSido} For a $k$-hypergraph pair $H$ over $V:=V_1
\times \ldots \times V_k$, a probability space ${\cal P}$, and a
zero-one function  $f \in L_H({\cal P})$ we have
$$\int f^H \ge \|f\|_1^{-|V_1|\ldots
|V_k|}.$$
\end{lemma}
\begin{proof}
Consider the $k$-hypergraph pair $K=(\frac{1}{2},\frac{1}{2})$ over
$V$. Lemma~\ref{lem:complete} shows that $K$ is a norming hypergraph
pair. Since $f$ is a zero-one function, we have $f^H \ge f^K$, and
thus by Corollary~\ref{cor:normOrder}
$$\int f^H \ge \int f^K \ge \|f\|_1^{|K|} \ge \|f\|_1^{-|V_1|\ldots
|V_k|}.$$
\end{proof}

\begin{lemma}
\label{lem:GowersApprox} Let $f,g: \Omega^k \rightarrow \mathbb{C}$
be two measurable functions with respect to the probability space
$(\Omega,{\cal F},\mu)$. Let $H=(\alpha,0)$ be a hypergraph pair
such that $\range(\alpha) \subseteq \{0,1\}$. Then

$$\left|\int f^H - g^H\right| \le |H| \|f-g\|_{U_k} \max(\|f\|_\infty,\|g\|_\infty)^{|H|-1}.$$
\end{lemma}
\begin{proof}
Let us label the elements of $\supp(\alpha)$ as
$\omega_1,\ldots,\omega_{|H|}$. Then for $0 \le i \le |H|$ define
$H_i := \sum_{j=1}^i 1_{\omega_j},$ so that $H_0=(0,0)$ and
$H_{|H|}=H$. Now by telescoping and applying
Corollary~\ref{cor:Gowers} we have
\begin{eqnarray*}
\left|\int f^H - g^H\right| &\le& \sum_{i=1}^{|H|} \left|\int
f^{H-H_{i-1}} g^{H_{i-1}} - f^{H-H_{i}}g^{H_{i}}\right| =
\sum_{i=1}^{|H|} \left|\int
f^{H-H_{i}}g^{H_{i-1}}(f^{1_{\omega_i}}-g^{1_{\omega_i}})\right|=
\\&=&\sum_{i=1}^{|H|} \left|\int
f^{H-H_{i}}g^{H_{i-1}}(f-g)^{1_{\omega_i}})\right|
  \le \sum_{i=1}^{|H|} \|f-g\|_{U_k} \|f\|_\infty^{|H|-i}
\|g\|_\infty^{i-1} \le \\ &\le&  |H| \|f-g\|_{U_k}
\max(\|f\|_\infty,\|g\|_\infty)^{|H|-1}.
\end{eqnarray*}

\end{proof}

\subsection{Proofs of Theorems~\ref{thm:outRuling}~and~\ref{thm:spreaded} \label{sec:ProofStructure}}

\begin{proof}[Theorem~\ref{thm:outRuling}]
Suppose that $H$ is a semi-norming $k$-hypergraph pair over $V=V_1
\times \ldots \times V_k$. The fact that $H \cong \overline{H}$
follows from Proposition~\ref{prop:isometry} because trivially
$|H|=|\overline{H}|$ and $\|f\|_H=\|f\|_{\overline H}$.

Now let $\epsilon>0$ be sufficiently small, and $h:\Omega^k
\rightarrow \{-1,1\}$ be such that $\|h\|_{U_k} \le \epsilon$ and
$\int h=0$, where here $(\Omega,{\cal F},\mu)$ is a
\emph{probability} space. The existence of $h$ is guaranteed by
Lemma~\ref{lem:RandomFunc}.

First we show that it is either the case that  for every $\psi \in
\supp(\alpha) \cup \supp(\beta)$, $\alpha(\psi)=\beta(\psi)$ or  for
every $\psi \in \supp(\alpha) \cup \supp(\beta)$,  $\{\alpha(\psi),
\beta(\psi)\}=\{0,1\}$, and we will handle the existence of a
universal $s$ later. Suppose that this statement fails for some
$\psi$. Note that at least one of $\alpha(\psi)$ or $\beta(\psi)$ is
not equal to $0$. We will assume that $\alpha(\psi) > \beta(\psi)$,
and the proof of the case $\alpha(\psi) <\beta(\psi)$ will be
similar. Since it is not the case that
$\beta(\psi)=1-\alpha(\psi)=0$, denoting $H-1_\psi =
(\alpha',\beta')$ we have
\begin{equation}
\label{eq:psiInSupp} \psi \in \supp(\alpha') \cup \supp(\beta').
\end{equation}
For $p:=\alpha(\psi)-\beta(\psi) \ge 0$, define $g:=h^{1/p}$, and
$$f:=\left\{ \begin{array}{lcl}1 & \qquad & h=1 \\ 0 & & h=-1  \end{array} \right. .$$
Since $\int h=0$, we have $\int f = 1/2$ and
\begin{equation}
\label{eq:Lowerfg} \int f^{H-1_\psi} g^{1_\psi} = \int f^H \ge
2^{-|V_1| \ldots |V_k|},
\end{equation}
where the equality follows from (\ref{eq:psiInSupp}) and the
definition of $f$, and the inequality follows from
Lemma~\ref{lem:subSido}. Denote by $K$ the hypergraph pair obtained
from $H$ by setting $\alpha(\psi)=\beta(\psi)=0$, i.e. $K := H -
\alpha(\psi) 1_\psi -  \beta(\psi) \overline{1_\psi}$. Now since
$|g|=1$, applying Corollary~\ref{cor:Gowers}, we have
$$ \left|\int g^H \right| =\left|\int g^K g^{\alpha(\psi)1_\psi+\beta(\psi)\overline{1_\psi}} \right|=
  \left|\int g^K |g|^{\beta(\psi)1_\psi} g^{p 1_\psi}\right| = \left|\int g^K h^{1_\psi}\right| \le
\|h\|_{U_k} \le \epsilon,$$
which shows that
\begin{equation}
\label{eq:boundgH} \|f\|_H^{H-1} \|g\|_H \le \|f\|_H^{H-1}
\epsilon^{1/|H|}.
\end{equation}
For sufficiently small $\epsilon$, (\ref{eq:Lowerfg}) and
(\ref{eq:boundgH}) contradict Lemma~\ref{lem:FirstHolder}.

Next we will prove the existence of a universal $s$. So suppose that
$H=(\alpha,\beta)$ is semi-norming and $\alpha=\beta$. Let $s=\max
\{ \alpha(\omega)+\beta(\omega): \omega \in V \}$. We will show that
$\frac{1}{s}H$ is semi-norming, and then
Corollary~\ref{cor:MoreThanOne} implies that
$\alpha(\omega)+\beta(\omega) \in \{0,s\}$. Let $\psi$ be such that
$\alpha(\psi)+\beta(\psi)=s$, and let $\tilde{H}_\psi =
\frac{1_\psi+\overline{1_\psi}}{2}$. Consider a measure space ${\cal
M}=(\Omega, {\cal F},\mu)$ and measurable functions $f,g:\Omega^k
\rightarrow \mathbb{C}$, and note that
\begin{eqnarray*}
\left| \int f^{(\frac{1}{s}H)-1_\psi} g^{1_\psi} \right| &\le& \int
|f|^{(\frac{1}{s}H)-1_\psi} |g|^{1_\psi} = \int
\left(|f|^{1/s}\right)^{H-s\tilde{H}_\psi}
\left(|g|^{1/s}\right)^{s\tilde{H}_\psi} \le \\ &\le&
\||f|^{1/s}\|_{H}^{|H|-s} \||g|^{1/s}\|_{H}^{s}=
\|f\|_{\frac{1}{s}H}^{\frac{1}{s}|H|-1} \|g\|_{\frac{1}{s}H},
\end{eqnarray*}
where in the second inequality we used Lemma~\ref{lem:Holder}. Now
Lemma~\ref{lem:FirstHolder} shows that $\frac{1}{s}H$ is a
semi-norming hypergraph pair, and this finishes the proof.
\end{proof}

Next we give the proof of Theorem~\ref{thm:spreaded}.

\noindent
\begin{proof}[Theorem~\ref{thm:spreaded}]
Suppose that $H= \dot\cup_{i=1}^m H_i$ where $H_i$ are
non-factorizable. Define $f:[0,1]^k \rightarrow \mathbb{R}$ as in
the following: $f(x_1,\ldots,x_k)=1$ if $\lfloor kx_1 \rfloor=
\ldots =\lfloor kx_k \rfloor$, and $f(x)=0$ otherwise. Then by
Corollary~\ref{cor:normOrder} we have
\begin{equation}
\label{eq:Diagonal} \|f\|_{H'} \le \|f\|_H.
\end{equation}
It is easy to see that
$$\int f^{H'} \ge k \left(\frac{1}{k}\right)^{|W_1|+\ldots+|W_k|} \times m,$$
while
$$\int f^{H} = k \left(\frac{1}{k}\right)^{|V_1|+\ldots+|V_k|} \times m.$$
Plugging these into~(\ref{eq:Diagonal}), and simplifying it, we
obtain the assertion of the theorem.
\end{proof}

\section{Geometry of the Hypergraph Norms \label{sec:Geometry}}

\subsection{Moduli of Smoothness and Convexity} Let us start by
recalling the definition of moduli of smoothness and convexity of a
normed space. For a normed space $X$, define the modulus of
smoothness as the function
\begin{equation}
\rho_X(\tau) = \sup \left\{ \frac{\|x - \tau y\|+ \|x + \tau
y\|}{2}-1 : \|x\|=\|y\|=1\right\},
\end{equation}
and the modulus of convexity as
\begin{equation}
\delta_X(\epsilon) = \inf \left\{ 1- \left\|\frac{x+y}{2}\right\| :
\|x\|=\|y\|=1, \|x-y\| \ge 2\epsilon\right\},
\end{equation}
where $0 \le \epsilon \le 1$. It should be noticed that the function
$\delta_X$ is frequently defined with $\epsilon$ in place of
$2\epsilon$. The following observation of
Lindenstrauss~\cite{Lindenstrauss63} shows that these two functions
behave in a dual form via Legendre transform:
\begin{equation}
\rho_{X^*}(\tau) = \sup \left\{ \tau \epsilon - \delta_X(\epsilon) :
0 \le \epsilon \le 1 \right\},
\end{equation}
where $X^*$ is the dual of $X$.

A normed space $X$ is called \emph{uniformly smooth}, if $\lim_{\tau
\rightarrow 0} \rho_X(\tau)/\tau = 0$, and it is called
\emph{uniformly convex}, if for every $\epsilon>0$,
$\delta_X(\epsilon)>0$. For $t \in (1,2]$ a normed space $X$ is said
to be \emph{$t$-uniformly smooth}, if there exists a constant $C>0$
such that $\rho_X(\tau) \le (C \tau)^t$, and for $r \in [2,\infty)$,
a normed space is said to be \emph{$r$-uniformly convex}, if there
exists a constant $C>0$ such that $\delta_X(\epsilon) \ge
(\epsilon/C)^r$. It is known that
$\rho_{\ell_2}(\tau)=(1+\tau^2)^{1/2}-1=\tau^2/2+O(\tau^4)$,
$\tau>0$  and
$\delta_{\ell_2}(\epsilon)=1-(1-\epsilon^2)^{1/2}=\epsilon^2/2+O(\epsilon^4)$
for $0<\epsilon<1$. Dvoretzky's theorem (see for
example~\cite{MilmanSchechtman}) implies that for every infinite
dimensional normed space $X$, we have $\rho_X(\tau) \ge
\rho_{\ell_2}(\tau)$ and $\delta_X(\epsilon) \le
\delta_{\ell_2}(\epsilon)$, and this was the reason for requiring $t
\in (1,2]$ and $r \in [2,\infty)$ in the definition of $t$-uniform
smoothness and $r$-uniform convexity. We will give another
equivalent definition for the notions of $t$-uniform smoothness and
$r$-uniform convexity due to Ball {\it et al}~\cite{BallCarlenLieb}.
First we need two simple lemmas.
\begin{lemma}
\label{lem:BonamiBeckner} Let $1<p \le q < \infty$ and
$\rho=\sqrt{\frac{p-1}{q-1}}$. Then for every two vectors $x$ and
$y$ in an arbitrary normed space $X$, we have

$$\left(\frac{\|x+\rho y\|^q + \|x-\rho y\|^q}{2} \right)^{1/q} \le \left(\frac{\|x+y\|^p+\|x-y\|^p}{2} \right)^{1/p}.$$
\end{lemma}
For the proof of Lemma~\ref{lem:BonamiBeckner} see Corollary 1.e.14
in~\cite{ClassicalBanach}.
\begin{lemma}
\label{lem:hypercont} Let $t \in (1,2]$, $r \in [2,\infty)$, and
$1<p,q<\infty$. Then there exists constants  $C=C(t,p)$ and
$C^*=C^*(r,q)$ such that for every $x,y \in \mathbb{C}$,
\begin{equation}
\label{eq:hyper1} \left(\frac{|x+y|^{p}+|x-y|^{p}}{2} \right)^{1/p}
\le \left(|x|^{t} + |C y|^{t}\right)^{1/t},
\end{equation}
and
\begin{equation}
\label{eq:hyper2} \left(\frac{|x+y|^{q}+|x-y|^{q}}{2} \right)^{1/q}
\ge \left(|x|^{r} + \left|\frac{1}{C^*} y\right|^{r}\right)^{1/r}.
\end{equation}
Furthermore for the best constants one can assume $C(t,p)=C^*(r,q)$,
if $\frac{1}{r}+\frac{1}{t}=1$ and $\frac{1}{p}+\frac{1}{q}=1$.
\end{lemma}
\begin{proof}
We only prove (\ref{eq:hyper1}), and (\ref{eq:hyper2}) as well as
the last assertion of the lemma will follow from duality by
Proposition~\ref{prop:duality} below. It suffices to prove the
theorem for $t=2$ as the right-hand side of (\ref{eq:hyper1}) is a
decreasing function in $t$. By Lemma~\ref{lem:BonamiBeckner}, we
have
$$\left(\frac{|x+y|^{p}+|x-y|^{p}}{2} \right)^{1/p} \le \left(\frac{|x+\rho y|^{2}+|x-\rho y|^{2}}{2} \right)^{1/2}
\le (|x|^2+|\rho y|^2)^{1/2},$$
where $\rho = \max(1,\sqrt{p-1})$.
\end{proof}

Now for a normed space $X$, inspired by Lemma~\ref{lem:hypercont},
for $1 < t \le 2 \le r < \infty$, and $1 < p,q < \infty$, one can
investigate the validity of the following two inequalities:
\begin{equation}
\label{eq:t,p-smoothness} \left(\frac{\|x+y\|^{p}+\|x-y\|^{p}}{2}
\right)^{1/p} \le \left(\|x\|^{t} + \|K y\|^{t}\right)^{1/t},
\end{equation}
and
\begin{equation}
\label{eq:r,q-convexity} \left(\frac{\|x+y\|^{q}+\|x-y\|^{q}}{2}
\right)^{1/q} \ge \left(\|x\|^{r} + \|K^{-1} y\|^{r}\right)^{1/r}
\end{equation}
where $K$ is a constant. We denote the smallest constant $K$ such
that (\ref{eq:t,p-smoothness}) is satisfied for all $x,y \in X$ by
$K_{t,p}(X)$ and similarly the smallest constant such that
(\ref{eq:r,q-convexity}) is satisfied by $K^*_{r,q}(X)$. Trivially
$K_{t,p}(X) \ge C(t,p)$ and $K^*_{r,q}(X) \ge C^*(r,q)$ where
$C(t,p)$ and $C^*(r,q)$ are the constants defined in
Lemma~\ref{lem:hypercont}.

\begin{remark}
\label{rem:constants} In the sequel $C(t,p)$ and $C^*(r,q)$ always
refer to the constants from Lemma~\ref{lem:hypercont}. Note that
$C(t,p)$ and $K_{t,p}(X)$ are both increasing in $t$ and $p$, and
$C^*(r,q)$ and $K^*_{r,q}(X)$ are both decreasing in $r$ and $q$.
Since Lemma~\ref{lem:BonamiBeckner} is valid for every normed space
$X$, for $1<p_2 \le p_1< \infty$,
\begin{eqnarray*}
\left(\frac{\|x+y\|^{p_1}+\|x-y\|^{p_1}}{2} \right)^{1/p_1} &\le&
\left(\frac{\left\|x+\sqrt{\frac{p_1-1}{p_2-1}}y\right\|^{p_2}+\left\|x-\sqrt{\frac{p_1-1}{p_2-1}}y\right\|^{p_2}}{2}
\right)^{1/p_2} \\ &\le& \left(\|x\|^{t} + \left\|K_{t,p_2}(X)
\sqrt{\frac{p_1-1}{p_2-1}} y\right\|^{t}\right)^{1/t},
\end{eqnarray*}
which implies $K_{t,p_1}(X) \le \sqrt{\frac{p_1-1}{p_2-1}}
K_{t,p_2}(X)$. Similarly for $1 < q_2 \le q_1<\infty$,
\begin{eqnarray*}
\left(\frac{\|x+y\|^{q_2}+\|x-y\|^{q_2}}{2} \right)^{1/q_2} &\ge&
\left(\frac{\left\|x+\sqrt{\frac{q_2-1}{q_1-1}}y\right\|^{q_1}+\left\|x-\sqrt{\frac{q_2-1}{q_1-1}}y\right\|^{q_1}}{2}
\right)^{1/q_1} \\ &\ge& \left(\|x\|^{r} +
\left\|\frac{1}{K^*_{r,q_1}(X)} \sqrt{\frac{q_2-1}{q_1-1}}
y\right\|^{r}\right)^{1/r},
\end{eqnarray*}
which shows that $K^*_{r,q_2}(X) \le \sqrt{\frac{q_1-1}{q_2-1}}
K^*_{r,q_1}(X)$.
\end{remark}

The following proposition which follows from
Remark~\ref{rem:constants}, and Proposition 7
in~\cite{BallCarlenLieb} shows that one can use
(\ref{eq:t,p-smoothness}) and (\ref{eq:r,q-convexity}) to give an
alternative definition of $t$-uniform smoothness and $r$-uniform
convexity.
\begin{proposition}
\label{prop:SmoothnessConvexity} Let $X$ be a $t$-uniformly smooth
normed space. Then for every $1 <p < \infty$, we have $K_{t,p}(X) <
\infty$. Conversely if $K_{t,p}(X)<\infty$ for some $1 < p <
\infty$, then $X$ is $t$-uniformly smooth.

Similarly let $Y$ be an $r$-uniformly convex normed space. Then for
every $1 <q< \infty$, we have $K^*_{r,q}(Y) < \infty$. Conversely if
$K^*_{r,q}(Y)<\infty$ for some $1 < q < \infty$, then $Y$ is
$r$-uniformly convex.
\end{proposition}
The constants $K_{t,p}$ and $K^*_{r,q}$ behave nicely with respect
to the duality. The proof of the following proposition is identical
to the proof of Lemma~5 from~\cite{BallCarlenLieb}, and thus we omit
it.
\begin{proposition}
\label{prop:duality} Consider a normed space $X$ and its dual $X^*$.
Suppose that $\frac{1}{p}+\frac{1}{q}=1$ and
$\frac{1}{r}+\frac{1}{t}=1$. Then $K_{r,p}(X) = K^*_{t,q}(X^*)$.
\end{proposition}

The notion of uniform convexity is first defined by Clarkson
in~\cite{Clarkson}, where he studied the smoothness and convexity of
$L_p$ spaces. To this end he established four inequalities known as
Clarkson inequalities. Let $1 < p \le 2 \le q < \infty$ and
$\frac{1}{p}+\frac{1}{q}=1$. In our notation the Clarkson
inequalities are the following: $K_{p,p}(L_p)=1$,
$K^*_{q,q}(L_q)=1$, $K^*_{q,p}(L_p)=1$, and $K_{p,q}(L_q)=1$. The
first two are easier to prove and known as ``easy'' Clarkson
inequalities, and the latter two are known as ``strong'' Clarkson
inequalities. The following observation shows that the strong
Clarkson inequalities imply the easy Clarkson inequalities.

\begin{observation}
\label{obs:Clarkson} Let $1 < t \le 2 \le r < \infty$ be such that
$\frac{1}{t}+\frac{1}{r}=1$. Then $K_{t,r}(X)=1$ if and only if
$K^*_{r,t}(X)=1$.
\end{observation}
\begin{proof}
Suppose that $K_{t,r}(X)=1$. Then for every $x,y \in X$, we have
$$\left(\frac{\|x+y\|^r+ \|x-y\|^r}{2}\right)^{1/r} \le \left(\|x\|^t+\|y\|^t\right)^{1/t}.$$
Now consider $x',y' \in X$. Replacing $x$ and $y$ in the above
inequality, respectively with $\frac{x'+y'}{2}$ and
$\frac{x'-y'}{2}$ we get
$$\left(\frac{\|x'\|^r+ \|y'\|^r}{2}\right)^{1/r} \le \left(\left\|\frac{x'+y'}{2}\right\|^t+\left\|\frac{x'-y'}{2}\right\|^t\right)^{1/t}.$$
which simplifies to
$$\left(\|x'\|^r+ \|y'\|^r\right)^{1/r} \le \left(\frac{\|x'+y'\|^t+\|x'-y'\|^t}{2}\right)^{1/t},$$
showing that $K^*_{r,t}(X)=1$. The proof of the converse direction
is similar.
\end{proof}

Consider $1 < p \le 2 \le q <\infty$. As we have already seen in
Proposition~\ref{prop:SmoothnessConvexity}, Clarkson's inequalities
imply that $L_p$ and $L_q$ spaces are both $p$-uniformly smooth and
$q$-uniformly convex. However this is not in general the best
possible. The actual situation is the following. The $L_p$ spaces
are $p$-uniformly smooth and $2$-uniformly convex, and the $L_q$
spaces are $2$-uniformly smooth and $q$-uniformly convex. These
facts are proved by Hanner~\cite{Hanner} through the so called
Hanner inequality. For $1 < p \le 2$, we say that a normed space
satisfies the $p$-Hanner inequality, if
$$\|x+y\|^p + \|x-y\|^p \ge (\|x\|+\|y\|)^p + \left| \|x\| - \|y\| \right|^p,$$
and for $2 \le q < \infty$, it satisfies the $q$-Hanner inequality
if
$$\|x+y\|^q + \|x-y\|^q \le (\|x\|+\|y\|)^q + \left| \|x\| - \|y\| \right|^q.$$
It is shown in~\cite{BallCarlenLieb} that if $X$ satisfies the
$p$-Hanner inequality, then $X^*$ satisfies the $q$-Hanner
inequality where $\frac{1}{p}+\frac{1}{q}=1$. The following
proposition reveals the relation between the Hanner inequality and
the notions of uniform smoothness and uniform convexity.

\begin{proposition}
\label{prop:HannerConstants} If a normed space $X$ satisfies the
$t$-Hanner inequality for $1< t \le 2$, then for every $2 \le q <
\infty$, we have $K^*_{q,t}(X)= C^*(q,t)$, and for every $1 < p \le
t'$, we have $K_{t,p}(X)=1$ where $\frac{1}{t}+\frac{1}{t'}=1$.

Similarly if a normed space $X$ satisfies the $r$-Hanner inequality
for $2 \le r  < \infty$, then for every $1 < p \le 2$, we have
$K_{p,r}(X) = C(p,r)$, and for every $r' \le q < \infty$, we have
$K^*_{r,q}(X) = 1$, where $\frac{1}{r}+\frac{1}{r'}=1$.
\end{proposition}
\begin{proof}
Suppose that $X$ satisfies the $t$-Hanner inequality for $1 < t \le
2$. Consider $2 \le q < \infty$, and $x,y \in X$. By the $t$-Hanner
inequality
$$\left(\frac{\|x+y\|^t+ \|x-y\|^t}{2}\right)^{1/t} \ge
\left(\frac{(\|x\|+\|y\|)^t+ |\|x\|-\|y\||^t}{2}\right)^{1/t} \ge
\left(\|x\|^q+\left\|\frac{1}{C^*(q,t)} y\right\|^q\right)^{1/q},$$
which shows that $K^*_{q,t}(X) \le C^*(q,t)$. But from this, and
Observation~\ref{obs:Clarkson} we also get $K_{t,t'}(X)=1$ as
$K^*_{t',t}(X) \le C^*(t',t)=1$. Hence  for $1 < p \le t'$ we have
$K_{t,p}(X)=1$. The second assertion follows from the first one by
duality.
\end{proof}

Inequalities (\ref{eq:t,p-smoothness}) and (\ref{eq:r,q-convexity})
are first appeared in~\cite{BallCarlenLieb}, where for $q \ge 2$,
the equalities
$K_{2,q}(\ell_q)=K_{2,q}(S_q)=K_{2,2}(\ell_q)=K_{2,2}(S_q)=\sqrt{q-1}$
are proved, where $S_q$ corresponds to the $q$-trace norm.

\begin{proposition}
\label{prop:LpConstants} For $1<t \le 2 \le r < \infty$,  $1<t_1 \le
2 \le r_1<\infty$, and $1 < p < \infty$, we have
\begin{equation}
\label{eq:lpK} K_{t_1,p}(\ell_r)= \left\{
\begin{array}{lcl}
C(t_1,r) & & p \le r \\
C(t_1,p) \le \cdot \le C(t_1,r) \sqrt{\frac{p-1}{r-1}} & \qquad & p \ge r \\
\end{array} \right.
\end{equation}
and
\begin{equation}
\label{eq:lpK*} K^*_{r_1,p}(\ell_r)= \left\{
\begin{array}{lcl}
C^*(r_1,t) & \qquad & p \ge t \\
C^*(r_1,p) \le \cdot \le C^*(r_1,t) \sqrt{\frac{t-1}{p-1}}  & & p
\le t
\end{array} \right.
\end{equation}
In particular $K_{2,p}(\ell_r)=\max(\sqrt{p-1},\sqrt{r-1})$, and
$K^*_{2,p}\left(\ell_t\right)=\max\left(\sqrt\frac{1}{p-1},\sqrt\frac{1}{r-1}\right)$.
\end{proposition}
\begin{proof}
It suffices to prove (\ref{eq:lpK}), and then (\ref{eq:lpK*}) will
follow from duality. Since $\ell_r$ satisfies the $r$-Hanner
inequality, by Proposition~\ref{prop:HannerConstants} we have
$K_{t_1,r}(\ell_r)=C(t_1,r)$. Then it follows from
Lemma~\ref{lem:BonamiBeckner} that for $p \ge r$, $K_{t_1,p}(\ell_r)
\le C(t_1,r) \sqrt{\frac{p-1}{r-1}}$. Furthermore since
$K_{t_1,p}(\ell_r)$ is increasing in $p$, we have $K_{t_1,p}(\ell_r)
\le C(t_1,r)$, for $p \le r$. It remains to show that
$K_{t_1,p}(\ell_r) \ge C(t_1,r)$ for $p \le r$. Consider two complex
numbers $a$ and $b$, and let $x,y \in \ell_r$ be as $x=(a,a)$ and
$y=(b,-b)$. Then since
$\|x+y\|_r=\|x-y\|_r=(|a+b|^r+|a-b|^r)^{1/r}$, plugging these two
vectors in
$$\left(\frac{\|x+y\|_r^{p}+\|x-y\|_r^{p}}{2}
\right)^{1/p} \le \left(\|x\|_r^{t_1} + \|K_{t_1,p}(\ell_r)
y\|_r^{t_1}\right)^{1/t_1},
$$
we get
$$\left(\frac{|a+b|^r+|a-b|^r}{2}\right)^{1/r} \le \left(|a|^{t_1}+|K_{t_1,p}(\ell_r)
b|^{t_1}\right)^{1/t_1},$$
which shows that $K_{t_1,p}(\ell_r) \ge C(t_1,r)$.
\end{proof}

Let $1 < t \le 2 \le r <\infty$ with $\frac{1}{t}+\frac{1}{r}=1$.
The spaces $\ell_t$ and $\ell_r$ are respectively $2$-uniformly
convex and $2$-uniformly smooth. Proposition~\ref{prop:LpConstants}
determines the optimum value of all corresponding constants. In
terms of the constants corresponding to $t$-uniformly smoothness of
$\ell_t$ and $r$-uniformly convexity of $\ell_r$, by
Remark~\ref{rem:constants} and Clarkson's inequalities we have
$$K_{t,p}(\ell_t)=\left\{
\begin{array}{lcl}
1 & \qquad & p \le r \\
C(t,p) \le \cdot \le \sqrt\frac{p-1}{r-1} & \qquad & p \ge r\\
\end{array} \right.
$$
and
$$K^*_{r,p}(\ell_t)=\left\{
\begin{array}{lcl}
1 & \qquad & p \ge t \\
C^*(r,p) \le \cdot \le \sqrt\frac{t-1}{p-1} & \qquad & p \le t\\
\end{array} \right.
$$

\begin{figure}[t]
\begin{center}
\begin{tabular}{|c|c|c|c|c|}

\cline{1-5}  & $1 < p \le t$  & $t \le p \le 2$ & $2 \le p \le r$ & $r \le p <\infty $ \\
\hline

$C(2,p)$ & 1 & 1 & $\sqrt{p-1}$ & $\sqrt{p-1}$  \\ \hline

$C^*(2,p)$ & $\sqrt{\frac{1}{p-1}}$ & $\sqrt{\frac{1}{p-1}}$ & 1 & 1  \\
\hline

$K_{2,p}(\ell_r)$ & $\sqrt{r-1}$ & $\sqrt{r-1}$ & $\sqrt{r-1}$ &
$\sqrt{p-1}$  \\ \hline

$K^*_{2,p}(\ell_t)$ & $\sqrt{\frac{1}{p-1}}$ & $\sqrt{\frac{1}{t-1}}$ & $\sqrt{\frac{1}{t-1}}$ & $\sqrt{\frac{1}{t-1}}$  \\
\hline

$C(t,p)$ & 1 & 1 & 1 & $\le \sqrt\frac{p-1}{r-1} $
\\ \hline

$K_{t,p}(\ell_t)$ & 1 & 1 & 1 & $\le \sqrt\frac{p-1}{r-1} $
\\ \hline

$C^*(r,p)$ & $\le \sqrt{\frac{t-1}{p-1}}$ & 1 & 1 &  1 \\ \hline

$K^*_{r,p}(\ell_r)$ & $\le \sqrt{\frac{t-1}{p-1}}$ & 1 & 1 &  1 \\
\hline

$K_{t_1,p}(\ell_r)$ & $C(t_1,r)$ & $C(t_1,r)$ & $C(t_1,r)$ & $\le
C(t_1,r) \sqrt{\frac{p-1}{r-1}}$  \\ \hline

$K^*_{r_1,p}(\ell_t)$ & $\le C^*(r_1,t) \sqrt{\frac{t-1}{p-1}}$ & $C^*(r_1,t)$ & $C^*(r_1,t)$ & $C^*(r_1,t)$  \\
\hline

\end{tabular}
\end{center}
\caption{\label{fig} Here $1 < t \le 2 \le r < \infty$ are such that
$\frac{1}{t}+\frac{1}{r}=1$, and $1 < t_1 \le 2 \le r_1 < \infty$
are arbitrary. }
\end{figure}

The moduli of smoothness and  convexity of a Banach space are only
isometric invariant, and they may change considerably under an
equivalent renorming. This leads to the definition of type and
cotype. A normed space is of \emph{type} $1 \le t \le 2$ if there
exists a constant $T_t$ such that for every integer $n \ge 0$, and
every set of vectors $x_1,\ldots,x_n$,
$$\E \left\|\sum_{i=1}^n \epsilon_i x_i \right\|^t \le T_t \left( \sum_{i=1}^n \|x_i\|^t \right)^{1/t},$$
where $\epsilon_i$ are independent Bernoulli random variables taking
values uniformly in $\{-1,1\}$. Similarly a normed space is said to
be of \emph{cotype} $2 \le r \le \infty$ if there exists a constant
$C_r$ such that  for every integer $n \ge 0$, and every set of
vectors $x_1,\ldots,x_n$,
$$ \left( \sum_{i=1}^n \|x_i\|^r \right)^{1/r} \le C_r \E \left\|\sum_{i=1}^n \epsilon_i x_i \right\|,$$
where in the case $r=\infty$ the left hand-side must be replaced by
$\max_{i=1}^n \|x_i\|$.

Trivially every normed space is of type $1$ and of cotype $\infty$.
If a normed space is of type $t_0$ and cotype $r_0$, then it is also
of type $t$ and cotype $r$ provided that $t \le t_0 \le 2 \le r_0
\le r$. Note that type and cotype do not change under an equivalent
norm. Figiel and Pisier~\cite{Figiel,FigielPisier} proved that
$t$-uniform smoothness implies type $t$, and $r$-uniform convexity
implies cotype $r$. The reverse is of course not true as for example
every finite dimensional space is of type and cotype $2$.

For $\lambda \ge 1$, a normed space $X$ is said to be
$\lambda$-\emph{finitely representable} in a normed space $Y$, if
for every finite dimensional subspace $E \subseteq X$, there exists
a linear map $T:E \rightarrow Y$ such that $\|T\|\|T^{-1}\| \le
\lambda$. If for every $\lambda>1$, $X$ is $\lambda$-finitely
representable in $Y$, then we simply say $X$ is \emph{finitely
representable} in $Y$.

It is well-known that infinite dimensional $L_p$ spaces are of type
$\min(p,2)$ and cotype $\max(2,p)$, and nothing better. Thus if
$\ell_p$ is $\lambda$-finitely representable in an space $X$ of type
$t$ and cotype $r$, then $t \le \min(2,p)$ and $r \ge \max(2,p)$. A
beautiful theorem due to Maurey and Pisier~\cite{MaureyPisier} says
that the converse is also true, i.e.  $\ell_p$ and $\ell_q$ are
finitely representable in $X$ where $p = \sup \{t : \mbox{$X$ is of
type $t$}\}$ and $q = \inf \{r : \mbox{$X$ is of cotype $r$}\}$.

Thus in order to study the type, cotype, modulus of smoothness, and
modulus of convexity of a normed space $X$, it is natural therefore
to first try to find the smallest $p \ge 1$ and largest $q$ that
$\ell_p$ and $\ell_q$ are finitely representable in $X$.

For a hypergraph pair $H$, define $\ell_H := L_H(\mathbb{N})$ where
$\mathbb{N}$ is endowed with the counting measure.

\begin{theorem}
\label{thm:embedding} If $H=(\alpha,\beta)$ is a non-factorizable
semi-norming hypergraph pair, then $\ell_{|H|}$ is a subspace of
$\ell_{H}$. Furthermore if $H$ is of Type I with parameter $s \le
2$, then $\ell_{s}$ is finitely representable in $\ell_{H}$.
\end{theorem}

The first part of the theorem which is trivial, shows that any
infinite dimensional $L_H$ space is not of any cotype $q <
\min(2,|H|)$. The second part which is more interesting and was
unknown to the author in~\cite{GraphNorms} shows that if $H$ is of
Type I with parameter $s <2$, then every infinite dimensional $L_H$
space is not of any type $p>s$. In particular in the case $s=1$, an
infinite dimensional $L_H$ space has no nontrivial type, and is not
uniformly smooth and convex. The next theorem shows that every such
space is of cotype $\min(2,|H|)$ which is the best possible by
Theorem~\ref{thm:embedding}.

\begin{theorem}
\label{thm:Cotype} Let $H$ be a non-factorizable semi-norming
hypergraph pair of Type I, then $\ell_H$  is of cotype
$\min(2,|H|)$.
\end{theorem}

In Theorem~\ref{thm:Cotype}, only the case $s=1$ is interesting to
us, as for $s>1$ we will prove something stronger in
Theorem~\ref{thm:ConvexitySmoothness}. The key to prove
Theorem~\ref{thm:Cotype} is the following observation. Consider a
non-factorizable semi-norming $k$-hypergraph pair
$H=(\alpha,\alpha)$ of Type I over $V:=V_1 \times \ldots \times
V_k$,  and functions $f_1, f_2,\ldots ,f_n \in \ell_H$. Then
$$
 \sum_{i=1}^n f_i^H = \sum_{i=1}^n \prod_{\omega
\in V} |f_i \circ \omega|^{2\alpha(\omega)} \le \prod_{\omega \in V}
\left(\sum_{i=1}^n |f_i \circ \omega|^{|H|} \right)^{1/|H|}=
\left(\sum_{i=1}^n |f_i|^{|H|}\right)^\frac{H}{|H|},
$$
where in the inequality above we used the classical H\"older
inequality. Hence~\footnote{Inequality~(\ref{eq:Concavity}) says
that $\ell_H$ is $|H|$-concave as a Banach lattice when $H$ is of
Type~I. For the definition of Banach lattice convexity and concavity
we refer the reader to~\cite{ClassicalBanach}.}
\begin{equation}
\label{eq:Concavity} \sum_{i=1}^n \|f_i\|_H^{|H|} \le \left\|
\left(\sum_{i=1}^n |f_i|^{|H|}\right)^{1/|H|}\right\|_H.
\end{equation}
We will also need the following
inequality~\footnote{Inequality~(\ref{eq:Convexity})  says that
$\ell_H$ is $s$-convex as a Banach lattice
(see~\cite{ClassicalBanach}).} in the sequel:
\begin{equation}
\label{eq:Convexity}\left\|\left(\sum_{i=1}^n
|f_i|^{s}\right)^\frac{1}{s} \right\|_H = \left(\int
\left(\sum_{i=1}^n |f_i|^{s}\right)^\frac{H}{s}\right)^{1/|H|} =
\left\|\sum_{i=1}^n |f_i|^s \right\|_{H/s}^{1/s} \le
\left(\sum_{i=1}^n \left\| |f_i|^s \right\|_{H/s} \right)^{1/s} =
\left(\sum_{i=1}^n \left\| f_i \right\|_{H}^s \right)^{1/s},
\end{equation}
where we used the fact that $H/s$ is also norming. Now we can state
the proof of Theorem~\ref{thm:Cotype}.

\noindent
\begin{proof}[Theorem~\ref{thm:Cotype}]
Consider  functions $f_1,\ldots,f_n \in \ell_H$, and let
$m:=\max(|H|,2)$. By applying Minkowski's inequality, Khintchine's
inequality, and then (\ref{eq:Concavity}), there exists a constant
$C$ such that
\begin{eqnarray*}
\E \left\|\sum_{i=1}^n \epsilon_i f_i \right\|_H &=& \E \left\|
\left|\sum_{i=1}^n \epsilon_i f_i \right| \right\|_H  \ge \left\| \E
\left|\sum_{i=1}^n \epsilon_i f_i \right| \right\|_H \ge C \left\|
\left( \sum_{i=1}^n  |f_i|^2\right)^{1/2} \right\|_H \\ &\ge& C
\left\| \left( \sum_{i=1}^n |f_i|^m\right)^{1/m} \right\|_H \ge C
\left(\sum_{i=1}^n \|f_i\|_H^m\right)^{1/m}.
\end{eqnarray*}
\end{proof}

Now let us turn to the other hypergraph pairs, i.e. the ones which
are not of Type I with parameter $1$. From
Theorem~\ref{thm:embedding}, in terms of the four parameters type,
cotype, modulus of smoothness, and of convexity, the following
theorem is the strongest statement one can hope to prove about them,
and in particular implies Theorem~\ref{thm:Cotype} for $H$ of Type I
with parameter $s>1$.

\begin{theorem}
\label{thm:ConvexitySmoothness} Let $H$ be a non-factorizable
semi-norming hypergraph pair such that $|H| \ge 2$.
\begin{itemize}
\item If $H$ is of Type II or Type I with parameter $s \ge 2$, then
$\ell_H$ is $2$-uniformly smooth and $|H|$-uniformly convex;

\item If $H$ is of Type I with parameter $1<s \le 2$, then
$\ell_H$ is $s$-uniformly smooth and $|H|$-uniformly convex.
\end{itemize}
\end{theorem}
\begin{remark}
If $1 < |H| < 2$, then it is easy to see by the previous results
that $\|\cdot\|_H$ corresponds to the $L_p$ norm where $p=|H|$, and
thus the Banach space properties of the norm are well-understood.
The case $|H|=1$ is also trivial.
\end{remark}
As it is discussed above, the notions of $t$-uniform smoothness and
$r$-uniform convexity can be further refined by looking at the
constants $K_{t,p}$ and $K^*_{r,q}$. In proving
Theorem~\ref{thm:ConvexitySmoothness} we will try to obtain the best
possible constants. This is treated and discussed in more details in
Section~\ref{sec:SmoothnessConvexityLH}. Next we prove
Theorems~\ref{thm:embedding}.

\subsection{Proof of Theorem~\ref{thm:embedding}}
Define $T:\ell_{|H|} \rightarrow \ell_H$ as $T: a \mapsto f_a$,
where for $a=\{a_i\}_{i \in \mathbb{N}}$, $f_a:\mathbb{N}^k
\rightarrow \mathbb{C}$ is defined as
$$f_a(i_1,\ldots,i_k) = \left\{
\begin{array}{lcl}
a_i & \qquad & i_1=i_2=\ldots=i_k=i  \\
0 & & \mbox{otherwise}
\end{array}\right.
$$
Since $H$ is non-factorizable, it is easy to see that $T$ is an
isometry.

Next we show that $\ell_s$ is finitely representable in $\ell_H$.
Since $L_H([0,1])$ is finitely representable in $\ell_H$, it
suffices to find a map $T:\ell_s([n]) \rightarrow L_H([0,1])$ with
$\|T\|\|T^{-1}\| \le 1+\epsilon$, for every $n \in \mathbb{N}$ and
every $\epsilon>0$. To this end  we find $f_1,\ldots,f_n:[0,1]^k
\rightarrow \mathbb{C}$, such that for every $x=(x_1,\ldots,x_n) \in
\ell_s([n])$ with $\|x\|_s=n^{1/s}$,
$$1-\epsilon/4 \le \left\|\sum_{i=1}^n x_i f_i \right\|_{H} \le 1+\epsilon/4,$$
and then the map $T:\ell_s([n]) \rightarrow L_H([0,1])$ defined by
$T: e_i \mapsto f_i$, for $i \in [n]$, satisfies $\|T\|\|T^{-1}\|
\le \frac{1+\epsilon/4}{1-\epsilon/4} \le 1+\epsilon$, for
$\epsilon<1$.
An argument similar to the proof of Lemma~\ref{lem:RandomFunc},
shows  that there exists $f_1,\ldots,f_n : [0,1]^k \rightarrow
\{0,1\}$ such that $\sum f_i=1$, and for every $i \in [n]$, $\int
f_i=\frac{1}{n}$ and
 $\|f_i-\frac{1}{n}\|_{U_k} \le \delta$. Note that since
$f_i$ are zero-one valued functions, $\sum_{i=1}^n f_i=1$ implies
that the supports of $f_i$ are pairwise disjoint. Then we have
$$\int \left(\sum_{i=1}^n x_i f_i\right)^{H} = \int \left(\sum_{i=1}^n
|x_i|^s
f_i\right)^{\tilde{H}},$$
where $\tilde{H}=(\frac{\alpha+\beta}{s},0)$. Furthermore if
$\|x\|_s=n^{1/s}$, then
$$\left\|\left(\sum_{i=1}^n |x_i|^s f_i \right)- 1\right\|_{U_k}= \left\|\sum_{i=1}^n \left(|x_i|^s f_i - \frac{|x_i|^s}{n}\right)\right\|_{U_k}
\le \sum_{i=1}^n |x_i|^s \left\| f_i - \frac{1}{n}\right\|_{U_k} \le
\delta \|x\|_s = \delta n^{1/s}.$$
Now by Lemma~\ref{lem:GowersApprox}
$$\left| \int \left(\sum_{i=1}^n |x_i|^s
f_i\right)^{\tilde{H}} - 1 \right| = \left| \int \left(\sum_{i=1}^n
|x_i|^s f_i\right)^{\tilde{H}} - 1^{\tilde{H}} \right| \le \delta
n^{1/s} |\tilde{H}| \max\left(\left\|\sum_{i=1}^n |x_i|^s
f_i\right\|_\infty, 1\right)^{|\tilde{H}|-1} \le \delta
n^{|\tilde{H}|} |\tilde{H}|.$$
Now taking $\delta$ sufficiently small finishes the proof.

\subsection{Complex Interpolation}

Let us recall the definition of the complex interpolation spaces.
Two topological vector spaces are called compatible, if there exists
a Hausdorff topological vector space containing both of these spaces
as subspaces. Consider two compatible normed space $X_0$ and $X_1$
and endow the space $X_0+X_1$ with the norm
$\|f\|_{X_0+X_1}=\inf_{f=f_0+f_1} (\|f_0\|_{X_0}+\|f_1\|_{X_1})$.
For  every $0 \le \theta \le 1$, one constructs the corresponding
complex interpolation space $[X_0,X_1]_\theta$, as in the following.

Let ${\cal F}(X_0,X_1)$ be the set of all analytic function $v:\{z:
0 \le \Re z \le 1\} \rightarrow X_0+X_1$ which are continuous and
bounded on the boundary, and moreover the function $t \rightarrow
v(j+it)$ $(j=0,1)$ are continuous functions from the real line into
$X_j$ which tend to zero as $|t| \rightarrow \infty$. We provide the
vector space ${\cal F}$ with a norm
$$\| v \|_{\cal F}:= \max \left\{\sup_{x \in \mathbb{R}} \|v(ix)\|_{X_0}, \sup_{x \in \mathbb{R}}\|v(1+ix)\|_{X_1}\right\}.$$
Then for every $0 \le \theta \le 1$, the complex interpolation space
of $X_0$ and $X_1$ is a normed space $X_0 \cap X_1 \subseteq
[X_0,X_1]_\theta \subseteq X_0+X_1$ defined as
$$[X_0,X_1]_\theta := \{f \in X_0+X_1: v(\theta)=f \exists v \in {\cal
F}(X_0,X_1)|\},$$
with the following norm:
$$\|f\|_\theta = \inf \left\{\|v\|_{\cal F}: f=v(\theta), v \in {\cal F}(X_0,X_1)\right\}.$$
The space $[X_0,X_1]_\theta$ has an interesting property. Consider
compatible pairs $X_0, X_1$ and $Y_0, Y_1$. Let $T:X_0+X_1
\rightarrow Y_0+Y_1$ be a bounded linear map. Then
(see~\cite{InterpolationBook}),
\begin{equation}
\label{eq:interpolation} \|T\|_{[X_0,X_1]_\theta \rightarrow
[Y_0,Y_1]_\theta} \le \|T\|_{X_0 \rightarrow Y_0}^{1-\theta}
\|T\|_{X_1 \rightarrow Y_1}^\theta.
\end{equation}

\begin{theorem}
\label{thm:interpolation} Let ${\cal M}=(\Omega,{\cal F},\mu)$ be a
measure space and $H$ be a norming hypergraph pair of Type I with
parameter $1$. Then for every $0 \le \theta \le 1$, and
$\frac{1}{p}=\frac{1-\theta}{p_0}+\frac{\theta}{p_1}$, where
$p_0,p_1 \ge 1$,
$$[L_{p_0H}({\cal M}),L_{p_1H}({\cal M})]_\theta = L_{pH}({\cal M}).$$
\end{theorem}
\begin{proof}
Let $f:\Omega^k \rightarrow \mathbb{C}$ be a measurable function
with $\|f\|_{pH}=1$. Define
$$v:\{z: 0 \le \Re z \le 1\} \rightarrow L_{p_0H}({\cal M}) + L_{p_1H}({\cal M})$$
by
$$v(z)=|f|^{p(\frac{1-z}{p_0}+\frac{z}{p_1})}.$$
Then $v(\theta)=|f|$ which shows that
$$\|f\|_\theta \le \max \left\{ \sup_{x \in
\mathbb{R}} \|v(ix)\|_{p_0H}, \sup_{x \in \mathbb{R}}
\|v(1+ix)\|_{p_1H} \right\}.$$
But note that
$$\|v(ix)\|_{p_0H} = \left(\int |v(ix)|^{p_0H}\right)^{1/|p_0H|} =  \left(\int \left(|f|^{p/p_0}\right)^{p_0H}\right)^{1/|p_0H|}
=\left(\int |f|^{pH}\right)^{1/|p_0H|}=1,$$
and similarly $\|v(1+ix)\|_{p_1H} \le 1$ which shows that
$\|f\|_\theta \le \|f\|_{pH}$.

Now for the other direction assume that $\|f\|_{\theta}=1$. Then for
every $\epsilon>0$, there exists $v_\epsilon$ such that
$f=v_\epsilon(\theta)$ and $\|v_\epsilon\|_{\cal F} \le 1+\epsilon$.
By H\"older's inequality,
$$\|f\|_{pH}^{|H|} = \sup \left\{ \int f^H g^H: \|g\|_{qH} \le 1 \right\},$$
where $1=\frac{1}{p}+\frac{1}{q}$. Fix $g:\Omega^k \rightarrow
\mathbb{C}$ with $\|g\|_{qH} \le 1$, and define $$u:\{z: 0 \le \Re z
\le 1\} \rightarrow L_{q_0H}({\cal M}) + L_{q_1H}({\cal M})$$ by
$$u(z) = |g|^{q(\frac{1-z}{q_0}+\frac{z}{q_1})},$$
where $\frac{1}{q_0}+\frac{1}{p_0}=1$ and
$\frac{1}{q_1}+\frac{1}{p_1}=1$. Let
$$F_\epsilon(z) = \int v_\epsilon(z)^H u(z)^H,$$
and notice that
$$|F_\epsilon(ix)|= \int v_\epsilon(ix)^H u(ix)^H \le \|v_\epsilon(ix)\|^{|H|}_{p_0H} \|u(ix)\|^{|H|}_{q_0H}
\le \|v_\epsilon\|^{|H|}_{\cal F} \times \|g^{q/q_0}\|^{|H|}_{q_0H}
\le (1+\epsilon)^{|H|}.$$
Similarly
$$|F_\epsilon(1+ix)|= \int v_\epsilon(1+ix)^H u(1+ix)^H \le \|v_\epsilon(1+ix)\|^{|H|}_{p_1H} \|u(1+ix)\|^{|H|}_{q_1H}
\le \|v_\epsilon\|^{|H|}_{\cal F} \times \|g^{q/q_1}\|^{|H|}_{q_1H}
\le (1+\epsilon)^{|H|}.$$
Then
$$\left| \int f^H g^H \right|= \left| F_\epsilon(\theta) \right| \le 1+\epsilon,$$
which by tending $\epsilon$ to zero leads to $\|f\|_{pH} \le 1$. We
conclude that $\|f\|_{pH}=\|f\|_\theta$.
\end{proof}

\subsection{Proof of Theorem~\ref{thm:ConvexitySmoothness} \label{sec:SmoothnessConvexityLH}}

In this section we give sharp bounds on the moduli of smoothness and
convexity of the norms defined by semi-norming hypergraph pairs.
This of course will prove Theorem~\ref{thm:ConvexitySmoothness}.

Consider a non-factorizable semi-norming hypergraph pair $H$, and an
infinite dimensional space $L_H$. Theorem~\ref{thm:embedding} shows
that $L_H$ contains $\ell_{|H|}$ as a subspace, and thus
$K_{t,p}(\ell_{|H|}) \le K_{t,p}(L_H)$ and $K^*_{r,q}(\ell_{|H|})
\le K^*_{r,q}(L_H)$, for $1 < t \le 2 \le r <\infty$ and $1 < p,q <
\infty$. Comparing Proposition~\ref{prop:HannerConstants} with
Figure~\ref{fig} shows that proving the $|H|$-Hanner inequality for
$L_H$ spaces, gives the optimal values of $K_{2,p}(L_H)$ and
$K^*_{|H|,|H|}(L_H)$, for every $p
> 1$.

\begin{theorem}[Hanner Inequality]
\label{thm:Hanner} Let $H$ be a non-factorizable semi-norming
hypergraph pair which is either of Type II, or of Type I with an
even integer parameter. Then for every $f,g \in \ell_H$, we have
$$\|f+g\|_H^{|H|} + \|f-g\|^{|H|} \le (\|f\|_H + \|g\|_H)^{|H|} + \left|\|f\|_H - \|g\|_H\right|^{|H|}.$$
\end{theorem}
\begin{proof}
Without loss of generality assume that $\|f\|_H \ge \|g\|_H$. Let
${\cal H}$ be the set of all pairs $(H_1,H_2)$ such that $H_1$ and
$H_2$ are hypergraph pairs taking only nonnegative integer values,
and furthermore $H_1+H_2=H$ and $|H_2|$ is an even integer. Then
\begin{eqnarray*}
\|f+g\|_H^{|H|} + \|f-g\|^{|H|} &=& \int (f+g)^H + (f-g)^H =
\sum_{(H_1,H_2) \in {\cal H}} \int f^{H_1} g^{H_2} \\&\le&
\sum_{(H_1,H_2) \in {\cal H}} \|f\|_H^{|H_1|} \|g\|_H^{|H_2|} =
(\|f\|_H+\|g\|_H)^{|H|}+\left(\|f\|_H - \|g\|_H\right)^{|H|},
\end{eqnarray*}
where in the inequality we used Lemma~\ref{lem:Holder}.
\end{proof}

Consider a norming hypergraph pair $H$ of Type I with parameter $s <
2$ and $|H| \ge 2$. Note that for every $2 \le q < \infty$, $\ell_s$
does not satisfy the $q$-Hanner inequality, as otherwise it would be
$2$-uniformly convex. Hence it follows from
Theorem~\ref{thm:embedding} that $\ell_H$ does not satisfy the
$q$-Hanner inequality for any $2 \le q < \infty$. However we
conjecture the following.

\begin{conjecture}
\label{conj:Hanner} Let $H=(\alpha,\beta)$ be a non-factorizable
semi-norming hypergraph pair of Type I with parameter $s \ge 2$.
Then every $L_H$ space satisfies the $|H|$-Hanner inequality.
\end{conjecture}
Since we could not establish the $|H|$-Hanner inequality for all
norming hypergraph pairs of Type I we have to treat some of them
separately. The next two lemmas which give the optimum bounds for
uniform smoothness and convexity constants of $\ell_H$ when $H$ is a
non-factorizable hypergraph pair of Type I with parameter $s \ge 2$
would have been followed from a positive answer to
Conjecture~\ref{conj:Hanner}.

\begin{lemma}[2-Smoothness]
\label{lem:Smoothness} Let $H=(\alpha,\beta)$ be a non-factorizable
semi-norming $k$-hypergraph pair with $|H| \ge 2$. If $H$ is of Type
II, or of Type I with parameter $s \ge 2$, then
$$K_{2,p}(\ell_H) = K_{2,p}(\ell_{|H|})= \left\{
\begin{array}{lcl}
\sqrt{|H|-1} & & p \le |H| \\
\sqrt{p-1} & \qquad & p \ge |H| \\
\end{array} \right.
$$

\end{lemma}
\begin{proof}
If suffices to prove $K_{2,|H|}(\ell_H) \le \sqrt{|H|-1}$, and the
rest will follow from Remark~\ref{rem:constants}. Suppose that $H$
is defined over $V:=V_1 \times \ldots \times V_k$. For $f,g \in
\ell_H$, we have to prove
\begin{equation}
\left(\frac{\|f+g\|_H^{|H|}+\|f-g\|_H^{|H|}}{2}\right)^{2/|H|} \le
\|f\|_H^2 +(|H|-1) \|g\|_H^2.
\end{equation}
Consider the counting measure on $\{-1,1\}$, and define the two
functions $\epsilon_1,\epsilon_2:\{-1,1\}^k \rightarrow \{-1,0,1\}$
as
$$\epsilon_1(x_1,\ldots,x_k)=\left\{\begin{array}{lcl}1 & \qquad & x_1=\ldots=x_k \\
0 & & \mbox{otherwise} \end{array}\right.,$$
and
$$\epsilon_2(x_1,\ldots,x_k)=\left\{\begin{array}{lcl}x_1 & \qquad & x_1=\ldots=x_k \\
0 & & \mbox{otherwise}\end{array} \right. .$$
Note that since $H$ is non-factorizable, for $x \in \{-1,1\}^{V_1}
\times \ldots \times \{-1,1\}^{V_1}$, we have
\begin{eqnarray}
\label{eq:epsilon1}
\epsilon_1^H(x) = \left\{\begin{array}{lcl}1 & \qquad & x=(1,\ldots,1) \\
0 & & \mbox{otherwise} \end{array}\right.,
\end{eqnarray}
and
\begin{eqnarray}
\label{eq:epsilon2}
\epsilon_2^H(x) = \left\{\begin{array}{lcl}\eta & \qquad & x=(\eta,\ldots,\eta) \\
0 & & \mbox{otherwise} \end{array}\right.,
\end{eqnarray}
Let $\tilde{f}=f \otimes \epsilon_1$ and $\tilde{g} = g \otimes
\epsilon_2$. From~(\ref{eq:epsilon1}) and~(\ref{eq:epsilon2}) it is
easy to see that
$$\int (\tilde{f}+\tilde{g})^H = \int (\tilde{f}-\tilde{g})^H = \int (f+g)^H+(f-g)^H,$$
and $\int \tilde{f}^H = 2\int f^H$ and $\int \tilde{g}^H = 2 \int
g^H$. Hence it suffices to prove
$$\left(\frac{\int (\tilde{f}+\tilde{g})^H}{2} \right)^{2/|H|} \ge  \left(\frac{\int \tilde{f}^H}{2}\right)^{2/|H|}+(|H|-1) \left(\frac{\int \tilde{g}^H}{2} \right)^{2/|H|}.$$
which simplifies to
\begin{equation}
\label{eq:simplifieTildes} \left(\int
(\tilde{f}+\tilde{g})^H\right)^{2/|H|} \ge \left(\int
\tilde{f}^H\right)^{2/|H|}+(|H|-1) \left(\int
\tilde{g}^H\right)^{2/|H|}.
\end{equation}
We will show that for $0 \le t \le 1$
\begin{equation}
\label{eq:simplifieTildesWitht} \left(\int
(\tilde{f}+t\tilde{g})^H\right)^{2/|H|} \ge \left(\int
\tilde{f}^H\right)^{2/|H|}+t^2(|H|-1) \left(\int
\tilde{g}^H\right)^{2/|H|}.
\end{equation}
Note that (\ref{eq:simplifieTildesWitht}) reduces to
(\ref{eq:simplifieTildes}) for $t=1$. Consider the functions
$L,R:[0,1] \rightarrow \mathbb{R}$, defined as
$$L(t)=\left(\int(\tilde{f}+t\tilde{g})^H\right),$$
and
$$R(t)=\left(\int\tilde{f}^H\right)^{2/|H|}+t^2(|H|-1) \left(\int
\tilde{g}^H\right)^{2/|H|}.$$
We have
$$\frac{d}{dt}L(t) = \int \sum_{\psi \in V} \alpha(\psi)(\tilde{f}+t\tilde{g})^{H-1_\psi} \tilde{g}^{1_\psi} + \beta(\psi)
(\tilde{f}+t\tilde{g})^{H-\overline{1_\psi}}
\tilde{g}^{\overline{1_\psi}}.$$
Then
$$\frac{d}{dt}L(t)^{2/|H|} = \frac{2}{|H|}\left(\int \sum_{\psi \in V} \alpha(\psi)(\tilde{f}+t\tilde{g})^{H-1_\psi} \tilde{g}^{1_\psi} + \beta(\psi)
(\tilde{f}+t\tilde{g})^{H-\overline{1_\psi}}
\tilde{g}^{\overline{1_\psi}}\right)L(t)^{\frac{2-|H|}{|H|}}.$$
We want to compute the second derivative. Denote ${\cal H}=\{1_\psi:
\psi \in V\} \cup \{\overline{1_\psi}: \psi \in V\}$, and define
$\gamma:{\cal H} \rightarrow \mathbb{R}$ by $\gamma: 1_\psi \mapsto
\alpha(\psi)$ and $\gamma: \overline{1_\psi} \mapsto \beta(\psi)$.
We have
\begin{eqnarray*}
\frac{d^2}{dt^2}L(t)^{2/|H|} &=& \frac{2}{|H|}\left(\int \sum_{H_1
\neq H_2 \in {\cal H}} \gamma(H_1)
\gamma(H_2)(\tilde{f}+t\tilde{g})^{H-H_1-H_2} \tilde{g}^{H_1+H_2}
\right.
\\ && \left. + \sum_{H_1 \in {\cal H}} \gamma(H_1)
(\gamma(H_1)-1)(\tilde{f}+t\tilde{g})^{H-2H_1}
\tilde{g}^{2H_1}\right)L(t)^{\frac{2-|H|}{|H|}} +
\\ &&+\left(\frac{d}{dt}L(t)\right)^2
\frac{2(2-|H|)}{|H|^2}L(t)^{\frac{2-2|H|}{|H|}}.
\end{eqnarray*}
Recalling the definition of $\tilde{f}$ and $\tilde{g}$, it is easy
to see that
$$L(0)^{2/|H|}=R(0),$$
and since $ \int \tilde{f}^{H-1_\psi} \tilde{g}^{1_\psi} = \int
f^{H-1_\psi} g^{1_\psi} - \int f^{H-1_\psi} g^{1_\psi} =0$ and $\int
\tilde{f}^{H-\overline{1_\psi}} \tilde{g}^{\overline{1_\psi}} = \int
f^{H-\overline{1_\psi}} g^{\overline{1_\psi}} - \int
f^{H-\overline{1_\psi}} g^{\overline{1_\psi}} =0$, we have
$$\left.\frac{d}{dt}L(t)^{2/|H|}\right|_{t=0}=\left. \frac{d}{dt}R(t)\right|_{t=0}=0.$$
Furthermore since $H$ is of Type II or of Type I with parameter $s
\ge 2$, by Lemma~\ref{lem:Holder}, we have
\begin{eqnarray}
\nonumber
\frac{d^2}{dt^2}L(t)^{2/|H|}|_{t=0}&=&\frac{2}{|H|}\left(\int
\sum_{H_1 \neq H_2 \in {\cal H}} \gamma(H_1)
\gamma(H_2)\tilde{f}^{H-H_1-H_2} \tilde{g}^{H_1+H_2} \right. +
\\ \nonumber & & \left. \qquad \qquad \sum_{H_1 \in {\cal H}} \gamma(H_1)(\gamma(H_1)-1)
\tilde{f}^{H-2H_1} \tilde{g}^{2H_1} \right)L(0)^{\frac{2-|H|}{|H|}}
\\ \nonumber &\le& \frac{2}{|H|}\left(\sum_{H_1 \neq H_2 \in {\cal H}} \gamma(H_1)
\gamma(H_2)+\sum_{H_1 \in {\cal H}} \gamma(H_1)(\gamma(H_1)-1)
\right)(\|\tilde{f}\|_H^{|H|-2}
\|\tilde{g}\|_H^{2})\|\tilde{f}\|_H^{2-|H|} \\ &=&
2(|H|-1)\|\tilde{g}\|_H^2= \frac{d^2}{dt^2}R(t)|_{t=0}.
\label{eq:SecondDerAtZero}
\end{eqnarray}
Now for every $0 \le t_0 \le 1$, one can replace $\tilde{f}$ with
$\tilde{f}+t_0\tilde{g}$ in (\ref{eq:SecondDerAtZero}) and obtain
that for every $0 \le t_0 \le 1$
\begin{equation}
\label{eq:SecondDer} \frac{d^2}{dt^2}L(t)^{2/|H|}|_{t=t_0} \le
\frac{d^2}{dt^2}R(t)|_{t=t_0}.
\end{equation}
We conclude (\ref{eq:simplifieTildesWitht}).
\end{proof}

Next we prove  Clarkson's inequalities for $\ell_H$ when $H$ is a
semi-norming hypergraph pair of Type II or of Type I with parameter
$s \ge 2$. As it is mentioned above this would follow from
Conjecture~\ref{conj:Hanner}.

\begin{lemma}[Clarkson's Inequalities]
\label{lem:Clarkson} Let $H$ be a non-factorizable semi-norming
hypergraph pair of Type II or  Type I with parameter $s \ge 2$ such
that $q:=|H| \ge 2$. Then
$$K_{p,q}(\ell_H) =K^*_{q,p}(\ell_H) =K^*_{q,q}(\ell_H)=1,$$
where  $\frac{1}{p}+\frac{1}{q}=1$.
\end{lemma}
\begin{proof}
Recall that always $K^*_{q,q} \le K^*_{q,p}$. Hence it suffices to
prove $K_{p,q}(\ell_H)=1$, as by Observation~\ref{obs:Clarkson} this
would imply $K^*_{q,p}(\ell_H)=1$. To this end, we need to show that
for $f,g \in \ell_H$, we have
\begin{equation}
\label{eq:strongClarkson} \left( \frac{\left\|f+g\right\|_{H}^{q} +
\left\|f-g\right\|_{H}^{q}}{2} \right)^{1/q} \le \left(
\|f\|_{H}^{p} + \|g\|_{H}^{p} \right)^{1/p},
\end{equation}
which is equivalent to
\begin{equation}
\label{eq:strongClarkson2} \left(
\left\|\frac{f+g}{2}\right\|_{H}^{q} +
\left\|\frac{f-g}{2}\right\|_{H}^{q} \right)^{1/q} \le \left(
\frac{\|f\|_{H}^{p} + \|g\|_{H}^{p}}{2} \right)^{1/p}.
\end{equation}
Proposition~\ref{prop:HannerConstants} shows that
(\ref{eq:strongClarkson2}) follows from the $|H|$-Hanner inequality.
Hence Theorem~\ref{thm:Hanner} implies (\ref{eq:strongClarkson2})
when $H$ is of Type~II or it is of Type~I with parameter $s$ where
$s$ is an even integer. Next assume that $H$ is of Type~I with
parameter $s \ge 2$.

For a real $1 \le t < \infty$, and a norming hypergraph pair $G$,
define the norm $L_{t}(\ell_G)$ on the set of pairs $(f,g)$ where
$f,g \in \ell_G$ as
$$\|(f,g)\|_{L_t(\ell_G)}:= \left(\|f\|_{G}^t+\|g\|_{G}^t\right)^{1/t}.$$
Consider the linear map $T:(f,g) \mapsto
(\frac{f+g}{2},\frac{f-g}{2})$. Then (\ref{eq:strongClarkson2}) says
that
\begin{equation}
\label{eq:strongClarksonReform} \|T\|_{L_{p}(\ell_H) \rightarrow
L_{q}(\ell_H)} \le 2^{-\frac{1}{p}}.
\end{equation}
We will prove this by interpolation. Let $\tilde{H}=\frac{1}{s}H$,
and $s_0$ and $s_1$ be two even integers satisfying $2 \le s_0\le s
\le s_1$, and $\theta$ be such that
$\frac{1}{s}=\frac{1-\theta}{s_0}+\frac{\theta}{s_1}$. Then
$\frac{1}{p}=\frac{1-\theta}{t_0}+\frac{\theta}{t_1}$, where
$\frac{1}{t_0}+\frac{1}{s_0|\tilde{H}|}=1$ and
$\frac{1}{t_1}+\frac{1}{s_1|\tilde{H}|}=1$.
Theorem~\ref{thm:interpolation} above, together with Theorem~5.1.2
from~\cite{InterpolationBook} imply that
$$\left[L_{s_0|\tilde{H}|}(\ell_{s_0\tilde{H}}),L_{s_1|\tilde{H}|}(\ell_{s_1\tilde{H}})\right]_\theta=
L_{s|\tilde{H}|}(\left[\ell_{s_0\tilde{H}},\ell_{s_1\tilde{H}}\right]_\theta)
=L_{q}(\ell_H),$$
and
$$\left[L_{t_0}(\ell_{s_0\tilde{H}}),L_{t_1}(\ell_{s_1\tilde{H}})\right]_\theta=L_{p}(\left[\ell_{s_0\tilde{H}},\ell_{s_1\tilde{H}}\right]_\theta)
=L_{p}(\ell_H).$$
Furthermore
$$\left(2^{-\frac{1}{t_0}}\right)^{1-\theta}\left(2^{-\frac{1}{t_1}}\right)^{\theta}=2^{-\frac{1}{p}}.$$
Now since we know that (\ref{eq:strongClarksonReform}) holds for
even values of $s \ge 2$, we have
$$ \|T\|_{L_{t_0}(\ell_{s_0\tilde{H}}) \rightarrow
L_{s_0|\tilde{H}|}(\ell_{s_0\tilde{H}})} \le 2^{-\frac{1}{t_0}},
$$
and
$$ \|T\|_{L_{t_1}(\ell_{s_1\tilde{H}}) \rightarrow
L_{s_1|\tilde{H}|}(\ell_{s_1\tilde{H}})} \le 2^{-\frac{1}{t_1}}.
$$
Then interpolation (\ref{eq:interpolation}), implies
(\ref{eq:strongClarksonReform}).
\end{proof}

Next Lemma determines the moduli of smoothness and convexity of
non-factorizable semi-norming hypergraph pairs of Type I with
parameter $1<s \le 2$.

\begin{lemma}
\label{lem:finalSmoothnessConvexity} Let $H$ be a non-factorizable
semi-norming hypergraph pair of Type I with parameter $s>1$ with
$|H| \ge 1$. Then $K_{s,|H|}(\ell_H) = C(s,|H|)$ and $K^*_{|H|,s}(X)
= C^*(|H|,s)$.
\end{lemma}
\begin{proof}
Let $C:=C(s,|H|)$ and $C^*:=C^*(|H|,s)$. Consider $f,g \in \ell_H$.
By (\ref{eq:Concavity}) and (\ref{eq:Convexity}) we have
\begin{eqnarray*}
\left(\frac{\|f+g\|_H^{|H|} + \|f-g\|_H^{|H|}}{2}\right)^{1/|H|}
&\le&
\left\|\left(\frac{|f+g|^{|H|}+|f-g|^{|H|}}{2}\right)^{1/{|H|}}\right\|_H
\\ & \le &
\left\|\left(|f|^{s}+|C g|^s\right)^{1/s}\right\|_H
\\ &\le& \left(\|f\|^s_H + \|C g\|^s_H\right)^{1/s},
\end{eqnarray*}
which shows that $K_{s,|H|}(\ell_H) \le C$. To prove
$K^*_{|H|,s}=C^*$, note that by (\ref{eq:Convexity}) and
(\ref{eq:Concavity}) we have
\begin{eqnarray*}
\left(\frac{\|f+g\|_H^{s} + \|f-g\|_H^{s}}{2}\right)^{1/s} &\ge&
\left\|\left(\frac{|f+g|^{s} + |f-g|^{s}}{2}\right)^{1/s}\right\|_H
\\
& \ge &
\left\|\left(|f|^{|H|}+\left|\frac{1}{C^*}g\right|^{|H|}\right)^{1/|H|}
\right\|_H
\\ & \ge & \left(\|f\|_H^{|H|}+\left\|\frac{1}{C^*}g\right\|_H^{|H|}
\right)^{1/|H|}.
\end{eqnarray*}
\end{proof}

\begin{remark}
Note that all results in Section~\ref{sec:SmoothnessConvexityLH} are
stated for \emph{non-factorizable} semi-norming hypergraph pairs.
Consider a semi-norming hypergraph pair $H=H_1 \dot\cup \ldots
\dot\cup H_m$, where $H_i$'s are non-factorizable. If $H$ is of Type
I, then by Theorem~\ref{thm:factor}, $\|\cdot\|_H=\|\cdot\|_{H_1}$,
and thus one can apply the results of
Section~\ref{sec:SmoothnessConvexityLH} to $H_1$ instead. However
some of our results do not cover the case where $H$ is factorizable
and of Type II.
\end{remark}

\section*{Acknowledgements}
The author wishes to thank Bal{\'a}zs Szegedy for many valuable
discussions. He also wishes to thank Mikl\'os Simonovits for
pointing his attention to~\cite{Erdos}.

\bibliographystyle{plain}
\bibliography{norm}
\end{document}